\documentclass[preprint,12pt]{elsarticle}

\usepackage{amssymb}
\usepackage{amsmath}
\usepackage{algorithm}
\usepackage{algorithmicx}
\usepackage{algpseudocode}
\usepackage{amsthm}
\usepackage{graphicx}
\newtheorem{lemma}{Lemma}
\usepackage[table,xcdraw]{xcolor}
\newtheorem{theorem}{Theorem}
\usepackage{hyperref}
\hypersetup{
    colorlinks=true,
    linkcolor=blue,
    filecolor=magenta,      
    urlcolor=cyan,
    pdftitle={Overleaf Example},
    pdfpagemode=FullScreen,
    }
\usepackage{caption}
\DeclareCaptionLabelFormat{bold}{\textbf{#1 #2}}
\floatname{algorithm}{Alg.}

\begin{document}

\begin{frontmatter}



\title{A Symmetric Multigrid-Preconditioned Krylov Subspace Solver for Stokes Equations}


\author[inst1]{Yutian Tao}

\affiliation[inst1]{organization={The computer sciences department, University of Wisconsin-Madison},
            addressline={1210 W Dayton St}, 
            city={Madison},
            postcode={53706}, 
            state={Wisconsin},
            country={United States}}

\author[inst1]{Eftychios Sifakis}


\begin{abstract}
Numerical solution of discrete PDEs corresponding to saddle point problems is highly relevant to physical systems such as Stokes flow. However, scaling up numerical solvers for such systems is often met with challenges in efficiency and convergence. Multigrid is an approach with excellent applicability to elliptic problems such as the Stokes equations, and can be a solution to such challenges of scalability and efficiency. The degree of success of such methods, however, is highly contingent on the design of key components of the multigrid scheme, including the hierarchy of discretizations, and the relaxation scheme used. Additionally, in many practical cases, it may be more effective to use a multigrid scheme as a preconditioner to an iterative Krylov subspace solver, as opposed to striving for maximum efficacy of the relaxation scheme in all foreseeable settings. In this paper, we propose an efficient symmetric multigrid preconditioner for the Stokes Equations on a staggered finite difference discretization. Our contribution is focused on crafting a preconditioner that (a) is symmetric indefinite, matching the property of the Stokes system itself, (b) is appropriate for preconditioning the SQMR iterative scheme \cite{freund1994new}, and (c) has the requisite symmetry properties to be used in this context. In addition, our design is efficient in terms of computational cost and facilitates scaling to large domains.

\end{abstract}




\end{frontmatter}


\section{Introduction}
\label{sec:intro}
Saddle-point problems arise in many fields such as fluid dynamics \cite{elman2002preconditioners, de2007two}, structure mechanics \cite{farhat1992saddle, franceschini2019block} and magnetohydrodynamics \cite{phillips2016block}. As modern computational platforms advance, the demand for solving large-scale problems becomes more pronounced. However, the performance potential of modern computing platforms is often hindered by existing algorithms, which may lack the necessary efficiency in design. Multigrid methods \cite{trottenberg2000multigrid} are designed specifically for their potential scalability in handling large-scale problems with the advantage of linear time and space complexity in principle. Although the convergence rates of multigrid methods are -- in the best-case scenario -- independent of problem sizes, the design and optimization of multigrid components, such as relaxation schemes, can significantly impact this property. The relaxation scheme, often referred to as a smoother in multigrid, can be designed according to various classical techniques, including the distributive smoother \cite{oosterlee2006multigrid}, Uzawa smoother \cite{elman1994inexact}, Braess-Sarazin smoother \cite{braess1997efficient} and Vanka smoother \cite{vanka1986block}. The distributive smoother transforms the original equations into a right-preconditioned system, aiming to improve properties such as conditioning or numerical structure that facilitates the applicability of the relaxation scheme. This has been studied for solving problems like the Stokes equations \cite{oosterlee2006multigrid, bacuta2011new, wang2013multigrid} and the Oseen equations \cite{chen2015multigrid}. The Uzawa smoother transforms indefinite systems into positive definite formulations by using Schur complement, with applications in solving the Stokes equations \cite{maitre1985fast} and poroelasticity equations \cite{luo2017uzawa}. The Braess-Sarazin smoother is a variant of the pressure correction steps in SIMPLE-type algorithms \cite{patankar1983calculation}. It also relies on the approximation of Schur complement and has been introduced for tackling challenges such as the Stokes equations \cite{braess1997efficient, he2018local} and magnetohydrodynamic equations \cite{adler2016monolithic}. In contrast to the previous three smoothers, the Vanka smoother focuses on solving local overlapping saddle-point problems and updating several local degrees of freedom collectively. While the Vanka smoother is highly effective in reducing local residuals, it comes with higher per-iteration costs compared to other smoothing methods. It has proven effective in the Stokes equations \cite{vanka1986block, saberi2022restricted} and poroelasticity equations \cite{franco2018multigrid}. In addition to the choice of smoother, the convergence of multigrid also depends on factors like grid-operators \cite{niestegge1990analysis} and discretization methods \cite{brandt2011multigrid}.

Considering the sensitivity of multigrid methods, they are often more effective when used as a preconditioner to an iterative Krylov subspace solver. The multigrid-preconditioned conjugate gradient (MGPCG) method \cite{tatebe1993multigrid} is the prototypical approach for solving positive definite systems such as the Poisson equation. However, for saddle-point problems, a more general Krylov subspace solver is required, due to the indefiniteness of the discretized equations. The generalized minimal residual (GMRes) method \cite{saad1986gmres} emerges as the most common choice in this context, promising convergence for any asymmetric indefinite system. The multigrid-preconditioned GMRes has been employed in addressing transport equations \cite{oliveira1998preconditioned}, advection-diffusion equations \cite{oosterlee1996use, ramage1999multigrid} and Navier-Stokes equations \cite{anselmann2023efficiency}. Despite the indefiniteness of saddle-point problems, symmetry is present as a property in many situations such as the Stokes equation, Oseen equations and Helmholtz equation. The minimal residual (Minres) \cite{paige1975solution} method, designed for symmetric indefinite systems, consumes both less memory and less computational time than GMRes per iteration while maintaining similar convergence rates. Thus, it is a more appropriate option in such symmetric scenarios. However, Minres has the limitation that it can only be paired with positive definite preconditioners, making it incompatible with multigrid preconditioning (which for a problem like ours would yield a symmetric indefinite preconditioner). For this reason, we turn our attention to the symmetric quasi-minimal residual (SQMR) method \cite{freund1994new}, a variant of MINRES with identical work and storage requirements. SQMR supports symmetric indefinite preconditioners, opening up opportunities for the use of multigrid as the preconditioning scheme. To the best of our knowledge, no multigrid-preconditioned SQMR has been developed for symmetric saddle-point problems.

Since SQMR requires a symmetric (indefinite) preconditioner, care must be taken to preserve the symmetry property if multigrid is used in this context. This requires a number of design choices, many of which are trivial to implement (e.g. making sure that the prolongation and restriction are adjoint operators), but also requires a more delicate design of the relaxation scheme, so to not hinder the symmetry property. Inspired by the steps needed to preserve symmetry within the Gauss-Seidel method, we design an effective symmetric distributive smoother and Vanka smoother in our work. The main idea behind the design of a symmetric smoother involves relaxation in a specific order, followed by relaxation in the exact reverse order. For distributive relaxation, this also needs to be paired with a combination of left- and right-preconditioning in order to preserve symmetry. Our primary focus is on the multigrid-preconditioned SQMR method, utilizing a staggered finite difference discretization of the Stokes equations as the model problem. The main contributions of this work can be summarized as follows:

\begin{itemize}
\item We propose a multigrid-preconditioned SQMR method for symmetric saddle-point problems within the context of the Stokes equations.
\item We design two multigrid smoothers, based on distributive relaxation and Vanka relaxation respectively, taking care to preserve the symmetry of the operators. We combine these relaxation operators selectively, to achieve a balance between convergence rates and computational time per iteration.
\item We compare the performance of our new multigrid-preconditioned solver with the classical multigrid method on both 2D and 3D benchmarks.
\end{itemize}

The remaining structure of the paper is as follows: In section \ref{sec:mg}, we introduce the finite difference discretization of the Stokes equation and explain how we apply multigrid. In section \ref{sec:smoother}, we discuss how we ensure the symmetry of the distributive smoother and Vanka smoother. In section \ref{sec: precond}, we propose our approach to use multigrid as a preconditioner for SQMR. In section \ref{sec: discretization}, we describe our discrete domain design used for setting boundary conditions. In section \ref{sec:exp}, we compare the performance of our multigrid-preconditioned solver with the classical multigrid method on both 2D and 3D benchmarks, and also benchmark against un-preconditioned approaches. In section \ref{sec:additional_exp}, we present additional experiments to understand how problems of different resolutions and conditioning, and multigrid cycle schemes can affect our results. Finally, we draw our conclusions in section \ref{sec:conclusion}.

\section{Multigrid for the Stokes equations}
\label{sec:mg}
In this work, we consider the Stokes equations as written below, which apply to both 2D and 3D scenarios:
\begin{align}
\begin{split}
    \label{eq:stokes}
    -\eta\Delta\mathbf{u} + \nabla p = \mathbf{f}\\
    -\nabla \cdot \mathbf{u} = \mathbf{0}
\end{split}
\end{align}
where $\mathbf{u} = (u, v)^T$ (2D) or $\mathbf{u} = (u, v, w)^T$ (3D) is the fluid velocity vector field, $p$ is the scalar pressure field, $\mathbf{f}$ is the body force vector field and $\eta$ is the fluid viscosity, together with suitable boundary conditions. In this work, we treat the viscosity parameter $\eta$ as spatially constant.

\begin{figure}[htb]
\includegraphics[width=\columnwidth]{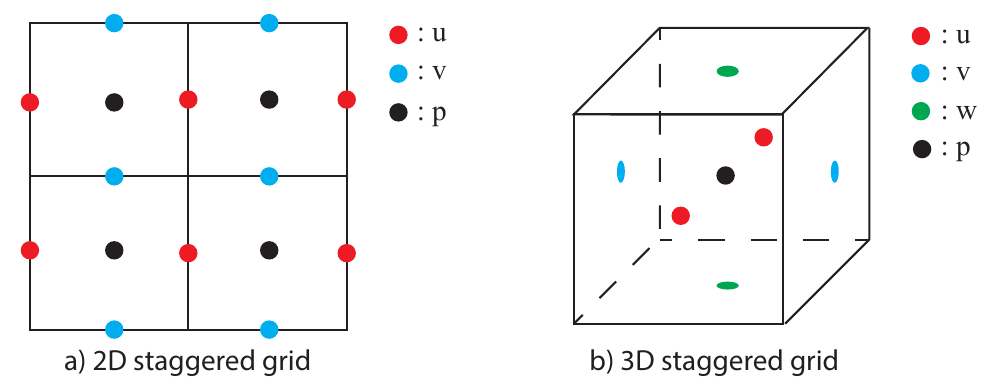}
\centering
\caption{The staggered grid discretization in 2D (left) and 3D (right)}
\label{fg: disc}
\end{figure}

For the discretization of the Stokes equations, we employ the standard marker-and-cell (MAC) staggered finite difference scheme \cite{trottenberg2000multigrid}. As shown in Fig. \ref{fg: disc}, the velocity components are centered at grid faces, while pressure variables are stored at cell centers. Although this discretization can be extended to non-uniform grid sizes $h_x \neq h_y \neq h_z$, we focus on uniform meshes with a grid size of $h=h_x=h_y=h_z$. The discretization employs central differences for the operators in the Stokes equations. Specifically, the Laplacian $\Delta$ operator is discretized using a five-point stencil in 2D and a seven-point stencil in 3D, while $\nabla$ operator is discretized via central difference approximation for velocity components. The $\nabla\cdot$ operator is discretized using central differences for pressure. This leads to the linear system discretized from the Stokes equations \ref{eq:stokes}:
\begin{equation}
    \mathbf{L}\mathbf{x} =
    \begin{pmatrix}
        \mathbf{A} & \mathbf{B}^T\\
        \mathbf{B} & \mathbf{0}\\
    \end{pmatrix}
    \begin{pmatrix}
        \mathbf{u}\\
        p\\
    \end{pmatrix}
    =
    \begin{pmatrix}
        \mathbf{f}\\
        0\\
    \end{pmatrix}
    = \mathbf{b}
\end{equation}
where $\mathbf{A}$, $\mathbf{B}$ and $\mathbf{B}^T$ represent discrete approximations of operators $-\eta\Delta$, $-\nabla\cdot$ and $\nabla$ operators respectively.

We proceed to describe the details and composition of our multigrid method. For our approach, which ultimately will be used to craft a preconditioner for an iterative Krylov solver, we adopt a V-cycle paradigm \cite{stuben2006multigrid}; while our contributions related to symmetric relaxation schemes are equally applicable to other types of multigrid cycles such as W-cycle.

\begin{figure}[htb]
\includegraphics[width=\columnwidth]{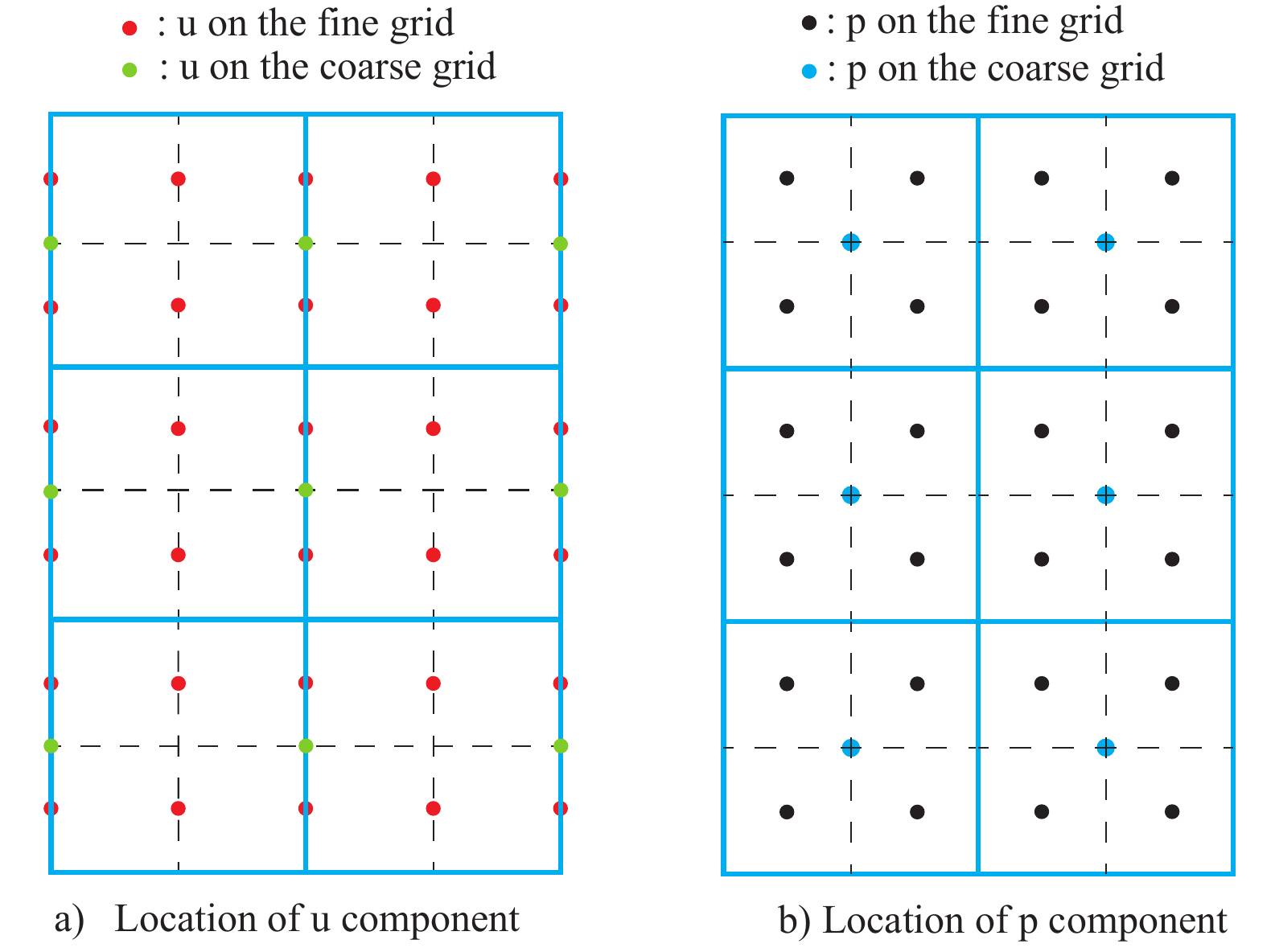}
\centering
\caption{2D discretization on fine and coarse level; just the horizontal component $u$ of the velocity field $\mathbf{u}=(u,v)$ is illustrated on the left, the placement of vertical velocity components $v$ is analogous, but on faces perpendicular to the $y$-axis.}
\label{fg: grid}
\end{figure}

The recursive V-cycle procedure is outlined in Alg. \ref{alg:vcycle}, where the superscript $h$ denotes the discretization grid size at every level of the multigrid hierarchy. We highlight that the operators $\mathbf{L}^{2h}, \mathbf{L}^{4h}, \ldots$ at coarser multigrid levels are constructed in our work based on a re-discretization of a coarsened description of the domain at successively coarser resolutions, as opposed to an algebraic Galerkin coarsening approach \cite{brandt2011multigrid} for example; for details of our domain coarsening approach see section \ref{sec: discretization}.
The locations of $u$ and $p$ components of the 2D fine coarse grids are shown in Fig. \ref{fg: grid}. The placement of variables and grids is analogous in 3D, with an 8-to-1 cell subdivision, and retaining face-centered velocities and cell-centered pressures.

\begin{figure}[htb]
\includegraphics[width=\columnwidth]{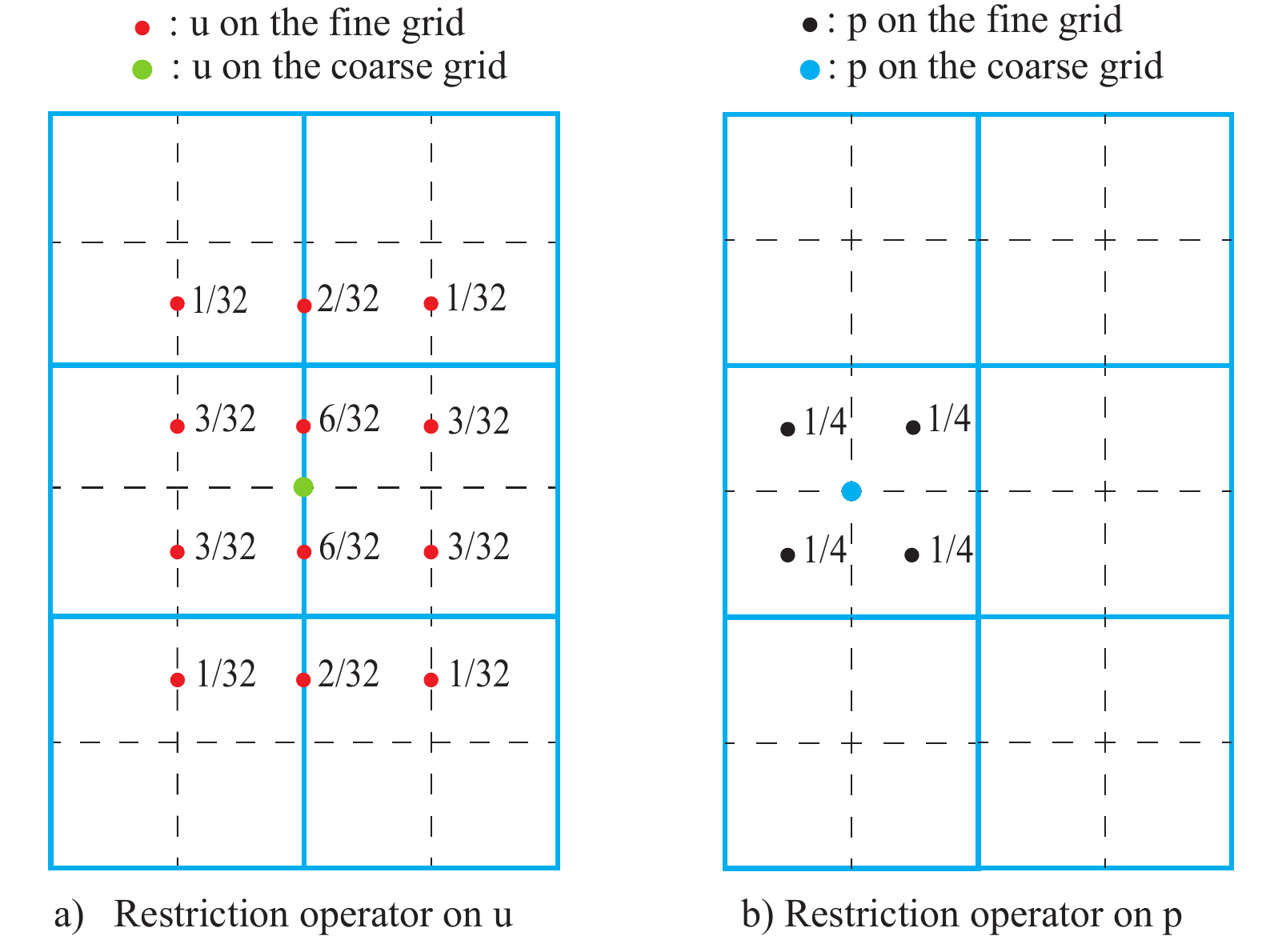}
\centering
\caption{2D restriction operator; just the horizontal component $u$ of the velocity field $\mathbf{u}=(u,v)$ is illustrated on the left, the placement of vertical velocity components $v$ is analogous, but on faces perpendicular to the $y$-axis.}
\label{fg: r_p}
\end{figure}

As we mentioned above, the restriction operator $\mathbf{R}$ and prolongation operator $\mathbf{P}$ can also influence convergence rates for the Stokes equations. Some common combinations have been investigated by Niestegge \cite{niestegge1990analysis}. 
Although this is not strictly necessary if multigrid is to be used as a solver outright, if we seek to use it as a \emph{symmetric} preconditioner, we must select restriction and prolongation operators $\mathbf{R}$ and $\mathbf{P}$ that are the transpose of each other up to scaling, meaning that $\mathbf{P}=c\mathbf{R}^T$ with $c$ typically being 4 for 2D and 8 for 3D. This scale factor (which does not impact the symmetry property) is typically selected to make both the prolongation and the restriction operators a partition of unity. To balance symmetry and the convergence rates, we use bilinear interpolation for velocities and constant per-cell interpolation for pressures as the prolongation operator. The restriction operator is thus defined as the transposed prolongation operator divided by $c$. Fig. \ref{fg: r_p} illustrates how the restriction stencil on $u$ component (the horizontal velocity component, the vertical one being analogous and symmetrical) and $p$ component in 2D cases.

\begin{algorithm}[htb]
\caption{Recursive V-cycle}
\begin{algorithmic}[1]
\State Input: System matrix $\mathbf{L}^h$, right hand side $\mathbf{b}^h$, initial guess $\mathbf{x}^h_0$ with grid size $h$ and multigrid level $n \geq 1$
\State Output: Solution $\mathbf{x}^h$
\Procedure{V-cycle} {$\mathbf{L}^h$, $\mathbf{b}^h$, $\mathbf{x}^h_0$, $n$}
  \If{$n = 1$}
    \State Solve $\mathbf{L}^h\mathbf{x}^h = \mathbf{b}^h$ exactly
  \Else
    \State $\mathbf{x}^h \gets \text{SMOOTH}(\mathbf{L}^h, \mathbf{x}^h_0, \mathbf{b}^h)$
    \State $\mathbf{r}^h \gets \mathbf{b}^h - \mathbf{L}^h\mathbf{x}^h$
    \State $\mathbf{r}^{2h} \gets \mathbf{R}^h\mathbf{r}^h$ \Comment{$\mathbf{R}^h$ is the restriction operator}
    \State $\mathbf{L}^{2h} \gets \text{Rediscretized from coarser level}$
    \State $\mathbf{e}^{2h} \gets \text{Recursively call } \text{V-CYCLE}(\mathbf{L}^{2h}, \mathbf{r}^{2h}, \mathbf{0}^{2h}, n - 1)$\;
    \State $\mathbf{e}^h \gets \mathbf{P}^h\mathbf{e}^{2h}$; \Comment{$\mathbf{P}^h$ is the prolongation operator}
    \State $\mathbf{x}^h \gets \mathbf{x}^h + \mathbf{e}^h$
    \State $\mathbf{x}^h \gets \text{SMOOTH}(\mathbf{L}^h, \mathbf{x}^h, \mathbf{b}^h)$\;
  \EndIf
\EndProcedure
\label{alg:vcycle}
\end{algorithmic}
\end{algorithm}

\section{Design of smoothing operator and preserving symmetry}
\label{sec:smoother}

Alg. \ref{alg:vcycle} details the procedure implementing a multigrid V-cycle. Iterative application of this cycle can be used to materialize an iterative solver, if used outright and not in the context of preconditioning. The same cycle can also be used to design a preconditioner appropriate for Krylov Subspace iterative methods. Tatebe \cite{tatebe1993multigrid} presented this process in the context of multigrid-preconditioned Conjugate Gradients, but the process by which the preconditioner is applied is more general across Krylov subspace solvers, and as follows: The preconditioner is intended to be a matrix $\mathbf{W}^\dag\approx\mathbf{L}^{-1}$ approximating the inverse of the system matrix $\mathbf{L}$. For our purposes the preconditioner must be a symmetric matrix, although we will use it in the context of a Krylov subspace solver that allows it to be indefinite; let us highlight that in the context of Stokes equations, an excellent preconditioner would easily be expected to be symmetric \emph{indefinite}, as the Stokes equations themselves are. 

As detailed in section \ref{sec: precond}, a matrix-vector multiplication $\mathbf{t}\leftarrow\mathbf{W}^\dag\mathbf{s}$ is implicitly implemented as 
$\mathbf{t}\leftarrow\textsc{V-cycle}(\mathbf{L}^h, \mathbf{b}^h\leftarrow \mathbf{s}, \mathbf{x}_0^h\leftarrow 0, n)$ via a call to the V-Cycle routine. Note that we have used the input vector $\mathbf{s}$ in the place of the right-hand-side vector $\mathbf{b}^h$, and used a zero initial guess $\mathbf{x}_0^h$ (the latter is essential in making sure the resulting operator is truly linear, as opposed to affine). Tatebe \cite{tatebe1993multigrid} identified the requirements for an operator $\mathbf{W}^\dag$ constructed in this way to be \emph{symmetric}. Two of the requirements are enforced by design, namely that the prolongation and restriction operators are transposed (up to scaling) of each other, and that the operator applied at the coarsest level of discretization is symmetric (here, by virtue of being an exact solution, i.e. an exact inverse of a symmetric matrix).

The more delicate requirement for the symmetry of $\mathbf{W}^\dag$ is that the smoothing operators at the restriction and prolongation phases of the V-cycle (lines 7 and 14 of Alg. \ref{alg:vcycle}) are transposes of each other. A straightforward way to enforce that would be to ensure that these smoothing operators are symmetric in their own right.
The smoothing operator can be formalized as:
\begin{equation}
    \label{eq:relaxation}
    \mathbf{x}^{k+1}=\mathbf{x}^{k}+\mathbf{S}(\mathbf{b}-\mathbf{L}\mathbf{x}^{k})
\end{equation}
where $k$ denotes the iteration step and $\mathbf{S}$ is the iteration matrix. For example, classical Gauss-Seidel smoothing can be expressed as $\mathbf{S}_{GS} = \mathbf{G}^{-1}$ with $\mathbf{G}$ being the lower triangular component of system matrix $\mathbf{L}$. Similarly, the weighted Jacobi method can be written as $\mathbf{S}_{WJ} = \omega\mathbf{D}^{-1}$, where $\omega$ is the damping factor and $\mathbf{D}$ is the diagonal component of $\mathbf{L}$.

Equation \ref{eq:relaxation} can be unrolled to yield $\mathbf{x}^k = \mathbf{C}\mathbf{b} + \mathbf{N}\mathbf{x}_0$. As detailed by Tatebe \cite{tatebe1993multigrid}, the symmetry of the smoothing operator for the purposes of crafting a symmetric preconditioner reduces to the symmetry of the matrix $\mathbf{C}$. We focus in our present work on the case where the same smoothing process is used in both the restriction and prolongation phases of the V-cycle; should it be an option to use different smoothing procedures at each phase, the corresponding $\mathbf{C}$ matrix for the two should be the transpose of each other; we do not focus on this case in our work though.

In the remainder of this section, we detail the steps that need to be taken such that this matrix $\mathbf{C}$ can be guaranteed to be symmetric, both when we use a distributive relaxation paradigm or a Vanka relaxation. 
We also explore how these operators can be combined while maintaining symmetry.

\subsection{Symmetric distributive smoother}
The design principle of a distributive smoother is that it can be viewed as a transformation of the original system $\mathbf{L}\mathbf{x}=\mathbf{b}$ into a right-preconditioned system $\mathbf{L}\mathbf{M}\mathbf{y}=\mathbf{b}$, where $\mathbf{M}\mathbf{y}=\mathbf{x}$. The hope is that the numerical properties of the combined matrix $\mathbf{LM}$ are more favorable for applying a standard relaxation procedure, rather than using the same procedure directly on $\mathbf{L}$. Furthermore, any relaxation on the system $\mathbf{LM}\mathbf{y}=\mathbf{b}$ can be implicitly emulated on the original system, without having to enact a change of variables from $\mathbf{x}$ to $\mathbf{y}$ beforehand. Specifically, any update $y_i \leftarrow y_i + \delta$ to one of the components of $\mathbf{y}$ can be emulated as $\mathbf{x}\leftarrow\mathbf{x}+\delta\mathbf{M}\mathbf{e}_i=\mathbf{x}+\delta\mathbf{m}_i$ on the $\mathbf{x}$, where $\mathbf{e}_i$ is the $i$-th basis vector (or, equally, $\mathbf{m}_i$ is the $i$-th column of $\mathbf{M}$)

For the Stokes equations, the transformed system can be formulated as
\begin{equation}
    \label{eq:dist}
    \mathbf{L}\mathbf{M}\mathbf{y} =
    \begin{pmatrix}
        \mathbf{A} & \mathbf{B}^T\\
        \mathbf{B} & \mathbf{0}\\
    \end{pmatrix}
    \begin{pmatrix}
        \mathbf{I} & -\mathbf{B}^T\\
        \mathbf{0} & \eta\mathbf{B}\mathbf{B}^T\\
    \end{pmatrix}
    \begin{pmatrix}
        \tilde{\mathbf{u}}\\
        \tilde{p}\\
    \end{pmatrix}
    = \mathbf{b}
\end{equation}
which yields
\begin{equation}
\begin{split}
    \label{eq:disc}
    \mathbf{L}\mathbf{M} &=
    \begin{pmatrix}
        \mathbf{A} & \mathbf{B}^T\\
        \mathbf{B} & \mathbf{0}\\
    \end{pmatrix}
    \begin{pmatrix}
        \mathbf{I} & -\mathbf{B}^T\\
        \mathbf{0} & \eta\mathbf{B}\mathbf{B}^T\\
    \end{pmatrix}
    =
    \begin{pmatrix}
        -\eta\Delta^h & \nabla^h\\
        -\nabla\cdot^h & \mathbf{0}\\
    \end{pmatrix}
    \begin{pmatrix}
        \mathbf{I} & -\nabla^h\\
        \mathbf{0} & -\eta\nabla\cdot^h\nabla^h\\
    \end{pmatrix}\\
    & =
    \begin{pmatrix}
        -\eta\Delta^h & \eta\Delta^h\nabla^h-\eta\nabla^h\Delta^h\\
        -\nabla\cdot^h & \Delta^h\\
    \end{pmatrix}
    =
    \begin{pmatrix}
        -\eta\Delta^h & \mathbf{0}\\
        -\nabla\cdot^h & \Delta^h\\
    \end{pmatrix}\\
\end{split}
\end{equation}
Here the superscript $h$ indicates the operator is discretized by a finite difference method with a grid size of $h$. The top right block is $\mathbf{0}$ due to two facts. First, $\nabla^h\cdot\nabla^h$ is equivalent to $\Delta^h$, which is a classical 5(2D) or 7(3D)-point Laplacian stencil in finite difference method. Secondly, the $\Delta$ and $\nabla$ operators are commutative, which is a property that Stokes equation specifically preserves under a finite difference discretization. It's important to note that when alternative discretization methods like Finite Elements are employed, the top-right block may not be exactly zero but only asymptotically zero up to truncation error. 
Furthermore, we should also note that the discrete matrix structure illustrated in equation \ref{eq:disc} is specific to the interior of the computational domain, since the presence of boundary equations could affect it (e.g. the triangular property of the combined matrix). 

The introduction of the distribution matrix $\mathbf{M}$ is motivated by the fact that $\mathbf{L}\mathbf{M}$ forms a block upper triangular matrix with the Laplacian operator on the diagonal block. This structure allows us to apply classical smoothing methods like Gauss-Seidel or weighted Jacobi on $\mathbf{L}\mathbf{M}$ matrix, leveraging their effectiveness on the Laplacian operator. The smoother for each component block will be applied in sequence. One way to appreciate this process would be as follows: If instead of a relaxation procedure, we applied a full solution for the Laplacian at each diagonal block (or, equivalently, if we smoothed to full convergence), iterating from one block component to the next would have been a block forward-substitution process. In practice, following the same relaxation schedule is a highly effective overall smoother, even if only a small number of relaxation steps are used at each block.

From that, $\delta\mathbf{y}$ is obtained to compute corrections $\delta\mathbf{x} = \mathbf{M}\delta\mathbf{y}$. The smoothing operator of, respectively,  distributive Gauss-Seidel (DGS) and distributive weighted Jacobi (DWJ) can be represented as $\mathbf{S}_{DGS} = \mathbf{M}\mathbf{\tilde{G}}^{-1}$ and $\mathbf{S}_{DWJ} = \omega\mathbf{M}\mathbf{\tilde{D}}^{-1}$ where $\mathbf{\tilde{G}}$ and $\mathbf{\tilde{D}}$ is the upper triangular and diagonal component of system matrix $\mathbf{L}\mathbf{M}$. Instead of computing the entire $\mathbf{L}\mathbf{M}$ matrix and applying Gauss-Seidel or weighted Jacobi, this process can be simplified by storing the diagonal part of $\mathbf{L}\mathbf{M}$ and using the connection $\mathbf{x}=\mathbf{M}\mathbf{y}$. Alg. \ref{alg:DGS} provides a detailed example of the simplification of DGS smoothing. A similar approach can be employed for DWJ smoothing.

\begin{algorithm}[htb]
\begin{algorithmic}[1]
\State Input: System matrix $\mathbf{L}$, initial guess $\mathbf{x}_0$ and right hand side $\mathbf{b}$
\State Output: Solution $\mathbf{x}$
\Procedure{DGS} {$\mathbf{L}$, $\mathbf{x}_0$, $\mathbf{b}$}
    \State Compute distributive matrix $\mathbf{M}$
    \State $\mathbf{d} = \mathbf{0}$
    \For{$j = 1, 2, ...$}
        \State $\mathbf{d}_j \gets \mathbf{L}_{j,} \cdot \mathbf{M}_{,j}$ 
        \State \Comment{where $\mathbf{L}_{j,}$ is the $j$th row of $\mathbf{L}$ and $\mathbf{M}_{,j}$ is the $j$th column of $\mathbf{M}$}
    \EndFor

    \State $\mathbf{x}^h \gets \mathbf{x}_0^h$
    \For{$j = 1, 2, ...$}
        \State $r \gets \mathbf{b}_j - \mathbf{L}_{j,}\mathbf{x}$
        \State $\mathbf{x} \gets \mathbf{x} + \frac{r}{\mathbf{d}_j}\mathbf{M}_{,j}$
    \EndFor
\EndProcedure
\caption{Distributive Gauss-Seidel Smoothing}
\label{alg:DGS}
\end{algorithmic}
\end{algorithm}

As mentioned above, we aim to design a symmetric smoother. However, the symmetry of distributive smoothing is hindered by right-preconditioning; we can easily show that applying an otherwise symmetric smoother (like Gauss-Seidel) to the non-symmetric distributed operator $\mathbf{LM}$ would compromise symmetry. In order to restore the symmetry property, we introduce left-preconditioning of the distributive smoother as follows
\begin{equation}
    \label{eq:ldist}
    \mathbf{M}^T\mathbf{L}\mathbf{x} =
    \begin{pmatrix}
        \mathbf{I} & \mathbf{0}\\
        -\mathbf{B} & \eta\mathbf{B}\mathbf{B}^T\\
    \end{pmatrix}
    \begin{pmatrix}
        \mathbf{A} & \mathbf{B}^T\\
        \mathbf{B} & \mathbf{0}\\
    \end{pmatrix}
    \begin{pmatrix}
        \mathbf{u}\\
        p\\
    \end{pmatrix}
    = \mathbf{M}^T\mathbf{b}
\end{equation}
It's obvious that $\mathbf{M}^T\mathbf{L} = (\mathbf{LM})^T$ due to the symmetry of $\mathbf{L}$. Our combined, symmetric distributive smoother operates as follows: In a first sweep, we choose a given traversal order of the degrees of freedom, and apply DGS according to Equation \ref{eq:dist} at each point. In a second sweep, we use exactly the reverse traversal order of the first step, and apply DGS, but in accordance with the left-preconditioned system in Equation \ref{eq:ldist}. The smoother constructed in this way is symmetric; a full proof of this fact is given in \ref{sec:pDGS}, and a similar proof applies to DWJ. Despite the efficiency of the distributive smoother, its convergence heavily relies on the commutativity of $\Delta$ and $\nabla$ operators in the interior region. This property does not discretely hold near the boundary, where this smoothing procedure is no longer applicable. Thus, we employ a different, more general and robustly convergent (albeit not as computationally efficient) smoother in the vicinity of the domain boundary.

\subsection{Symmetric Vanka smoother}
\begin{figure}[htb]
\includegraphics[width=7cm]{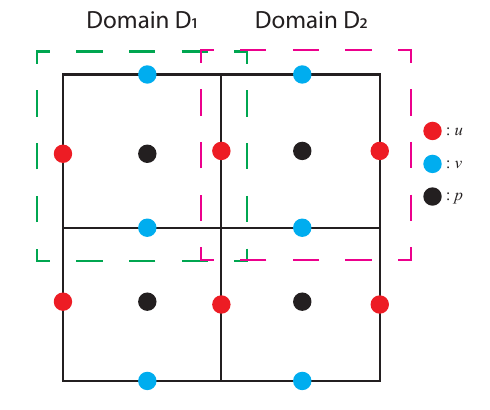}
\centering
\caption{Vanka Smoother on 2D example}
\label{fg: vanka}
\end{figure}

The Vanka smoother is more powerful and broadly applicable, but incurs a higher computational cost per iteration compared to the distributive smoother. It can be viewed as a specific type of Schwarz smoother \cite{schoberl2003schwarz}. The fundamental idea behind the Vanka smoother is to solve the Stokes equations locally, subdomain by subdomain, wherein all degrees of freedom within one subdomain are updated simultaneously. Like other smoothers, Vanka relaxation can be expressed in the format of Equation \ref{eq:relaxation}. In the Vanka smoother, we divide the entire discrete domain $D$ into $n$ potentially overlapping subdomains $\{D_1, ..., D_n\}$ as shown in Fig. \ref{fg: vanka}. Each subdomain can be viewed as a local saddle-point problem. Solving each subdomain problem sequentially, akin to Gauss-Seidel, is known as multiplicative Vanka relaxation. The smoother operator can be defined as
\begin{equation}
    \label{eq:mv}
    \mathbf{S}_{MV}=[\mathbf{I}-\prod_{i = 1}^{n}(\mathbf{I} - \mathbf{C}_i^T(\mathbf{C}_i\mathbf{L}\mathbf{C}_i^T)^{-1}\mathbf{C}_i\mathbf{L})]\mathbf{L}^{-1}
\end{equation}
where $\mathbf{C}_i$ represents the selection matrix from domain $D$ to domain $D_i$. Alternatively, solving each subdomain problem simultaneously as Jacobi method, is known as the additive Vanka smoother. The smoother operator for this approach is defined as
\begin{equation}
    \label{eq:av}
    \mathbf{S}_{AV}=\sum_{i=1}^{n}\mathbf{C}_i^T(\mathbf{C}_i\mathbf{L}\mathbf{C}_i^T)^{-1}\mathbf{C}_i
\end{equation}
We prove that the additive Vanka smoother is always symmetric in \ref{sec:pAV}. To ensure the symmetry in the multiplicative Vanka smoother, we adopt a similar approach to symmetric Gauss-Seidel. We solve each subdomain from domain $D_1$ to $D_n$, and then in reverse order from $D_n$ to $D_1$ again. The proof of symmetry for our symmetric multiplicative Vanka smoother is given in \ref{sec:pMV}.

\subsection{Symmetric integration of smoother}
As previously discussed, these two smoothers provide advantages depending on the region of the computational domain. The symmetric distributive smoother provides reliable convergence and attractive computational cost when applied to the interior region, while the symmetric Vanka smoother exhibits greater suitability near the boundary, albeit at a higher cost. However, despite the fact that each can be made symmetric on its own, it's not trivially true that any combination will inherently retain symmetry. To ensure strict symmetry with the multigrid framework, we propose the following smoothing step algorithm:

First, we partition the unknowns into two non-overlapping sets: boundary set $V_1$ and interior set $V_2$. Then we apply symmetric Vanka smoother to $V_1$, symmetric distributive smoother to $V_2$ and another symmetric Vanka smoother to $V_1$ again. The detailed algorithm is shown in Alg. \ref{alg:smoothing}, with a formal proof for its symmetry provided in \ref{sec:pVDV}. In practice for our implementation, DGS smoother is used instead of DWJ, because it achieves a more attractive convergence rate but a similar runtime cost. 

This approach leverages the strengths of both smoothers. On one hand, the larger size of $V_2$ ensures that the symmetric Vanka smoother isn't overused, enhancing computational efficiency. On the other hand, employing symmetric Vanka smoother near the boundary helps prevent the potential divergence issues of the symmetric distributive smoother. Therefore, this integrated approach is a balanced and effective choice in practical applications compared to using either smoother alone.

\begin{algorithm}[htb]
\begin{algorithmic}[1]
\State Input: System matrix $\mathbf{L}^h$, initial guess $\mathbf{x}^h_0$ and right hand side $\mathbf{b}^h$
\State Output: Solution $\mathbf{x}^h$
\Procedure{Smooth} {$\mathbf{L}^h$, $\mathbf{x}^h_0$, $\mathbf{b}^h$}
    \State Partition $\mathbf{x}^h$ into two non-overlapping sets: boundary set $V_1$ and interior set $V_2$
    \State $\mathbf{x}^h \gets \text{Symmetric\_Vanka}(\mathbf{L}^h, \mathbf{x}^h_0, \mathbf{b}^h)$ on the boundary set $V_1$
    \State $\mathbf{x}^h \gets \text{Symmetric\_Distributive}(\mathbf{L}^h, \mathbf{x}^h_0, \mathbf{b}^h)$ on the interior set $V_2$
    \State $\mathbf{x}^h \gets \text{Symmetric\_Vanka}(\mathbf{L}^h, \mathbf{x}^h_0, \mathbf{b}^h)$ on the boundary set $V_1$
\EndProcedure
\caption{Smoothing}
\label{alg:smoothing}
\end{algorithmic}
\end{algorithm}

\section{Preconditioning}
\label{sec: precond}
Although the use of the Vanka smoother significantly improves convergence near the boundary, multigrid is more stable when used as the preconditioner for a Krylov subspace solver. SQMR is known for its lower computational time and storage requirements, and supports arbitrary symmetric preconditioners, including indefinite ones from a multigrid V-cycle. Therefore, we propose using our multigrid as the preconditioner for SQMR. The application of the multigrid preconditioner is outlined in Alg. \ref{alg:mg-sqmr}. It's worth mentioning that not just any multigrid V-cycle can be used as the preconditioner. This procedure must be equivalent to $\mathbf{t} \gets \mathbf{M}^{-1}\mathbf{s}_0$ where matrix $\mathbf{M}$ is symmetric. Fortunately, as we prove in the appendix, our multigrid formulation indeed satisfies this condition with the initial guess $\mathbf{x}_0=\mathbf{0}$.

Another potential benefit of using multigrid as the preconditioner is that it allows us to use different matrices for SQMR and V-cycle, which means we can use some approximation $\tilde{\mathbf{L}}$ instead of $\mathbf{L}$ in our multigrid V-cycle. For example, in order to prevent $\mathbf{L}$ from failing to be full rank when applying all Dirichlet boundary conditions, we modify it as
\begin{equation}
    \tilde{\mathbf{L}} =
    \begin{pmatrix}
        \mathbf{A} & \mathbf{B}^T\\
        \mathbf{B} & -\gamma\mathbf{I}\\
    \end{pmatrix}
\end{equation}
where $\gamma$ is an adequately small number. This modification corresponds to the classical penalty method. Since the modification is only applied as a preconditioner instead of as a replacement of the system matrix, the divergence-free property can be always exactly maintained.

\begin{algorithm}[htb]
\begin{algorithmic}[1]
\State Input: System matrix $\mathbf{L}$, initial guess $\mathbf{x}_0=\mathbf{0}$, right hand side $\mathbf{b}$ and multigrid level $n \geq 1$
\State Output: Solution $\mathbf{x}$
\Procedure{MG-SQMR} {$\mathbf{L}$, $\mathbf{x}_0$, $\mathbf{b}$}
    \State $\mathbf{s}_0 \gets \mathbf{b} - \mathbf{L}\mathbf{x}_0$, $\mathbf{t} \gets \textcolor{red}{\text{V-CYCLE}(\tilde{\mathbf{L}}, \mathbf{s}_0, \mathbf{0}, n)}$, $\mathbf{q}_0 \gets \mathbf{t}$, $\tau_0 \gets ||\mathbf{t}||_2$, $\nu_0 \gets 0$, $\rho_0 \gets \mathbf{s}_0^T\mathbf{q}_0$, $\mathbf{d} \gets \mathbf{0}$\;

    \For{$j = 1, 2, ...$}
        \State $\mathbf{t} \gets \mathbf{L}\mathbf{q}_{j - 1}$, $\sigma_{j - 1} \gets \mathbf{q}_{j - 1}^T\mathbf{t}$, $\alpha_{j - 1} \gets \frac{\rho_{j - 1}}{\sigma_{j - 1}}$, $\mathbf{s}_j \gets \mathbf{s}_{j - 1} - \alpha_{j - 1}\mathbf{t}$
        \State $\mathbf{t} \gets \textcolor{red}{\text{V-CYCLE}(\tilde{\mathbf{L}}, \mathbf{s}_j, \mathbf{0}, n)}$, $\nu_j \gets \frac{||\mathbf{t}||_2}{\tau_{j - 1}}$, $c_j \gets \frac{1}{\sqrt{1 + \nu_j^2}}$, $\tau_j \gets \tau_{j - 1}\nu_j c_j$
        \State $\mathbf{d}_j \gets c_j^2\nu_{j - 1}^2 \mathbf{d}_{j - 1} + c_j^2\alpha_{j - 1}\mathbf{q}_{j - 1}$, $\mathbf{x}_j \gets \mathbf{x}_{j - 1} + \mathbf{d}_j$
        \If{$\mathbf{x}_j$ has converged}
            \State {$\mathbf{x} \gets \mathbf{x}_j$, stop}
        \EndIf
        \State $\rho_j \gets \mathbf{s}_j^T\mathbf{t}$, $\beta_j \gets \frac{\rho_j}{\rho_{j - 1}}$, $\mathbf{q}_j \gets \mathbf{t} + \beta_j\mathbf{q}_{j - 1}$
    \EndFor
\EndProcedure
\caption{Multigrid-Preconditioned SQMR}
\label{alg:mg-sqmr}
\end{algorithmic}
\end{algorithm}

\section{Discrete Domain Design}
\label{sec: discretization}
In this section, we describe our design decisions for the discrete description of the computational domains we target. We highlight that, as a matter of scope, we focus on a first-order accurate discretization of the Stokes equations, at grid-cell precision. This is primarily a choice for simplicity of exposition; our solver architecture and specifically our use of Vanka relaxation near the boundary could have been combined with a higher order scheme to provide sub-cell resolution of boundary conditions if desired. In addition, even if the finest level of discretization in the multigrid hierarchy was discretized at higher order, it would still be an option to use first-order accurate discretizations at coarser levels, especially when the multigrid cycle is intended as a preconditioner.

\begin{figure}[htb]
\includegraphics[width=\columnwidth]{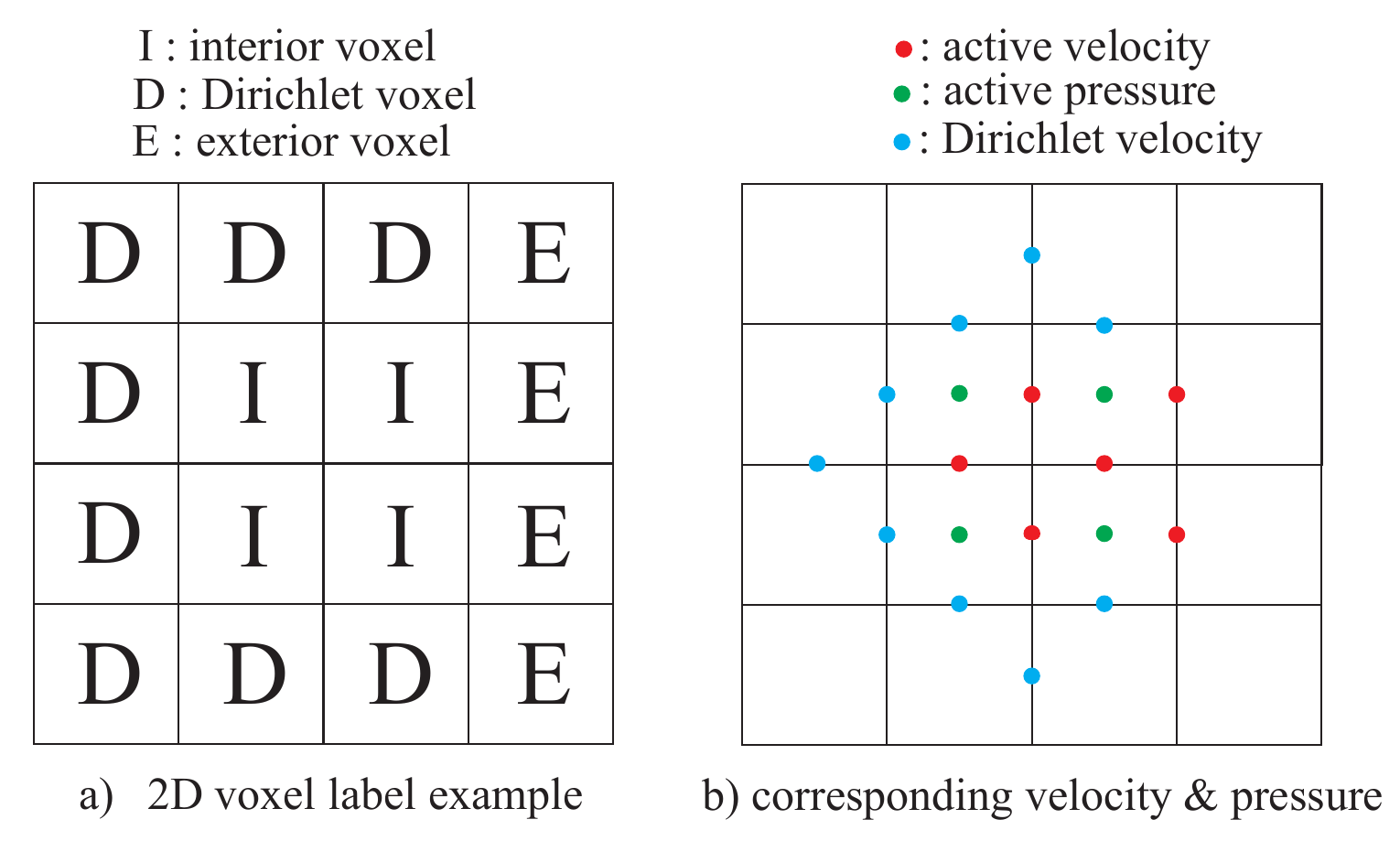}
\centering
\caption{2D Discrete Domain Example}
\label{fg: domain}
\end{figure}

We adopt a domain description where every cell is labeled as Dirichlet, exterior or interior as illustrated in Fig. \ref{fg: domain}. Using this designation of cells as input, we proceed to also classify velocity or pressure variables as ``Dirichlet'', ``active'' or ''inactive'' as well. Specifically, any velocity variables on faces of Dirichlet cells are treated as Dirichlet boundary conditions. Velocities on faces between two ``interior'' cells are treated as active degrees of freedom and are solved by our scheme. Any velocity variables that do not neighbor an active cell are treated as inactive, and excluded from our solve. Finally, pressures at centers of interior cells are active degrees of freedom, while others (at centers of Dirichlet or exterior cells) are inactive and excluded from the solve.

We adopt this cell-based designation, as opposed to, e.g. designating individual grid faces as Neumann faces, because our design leads itself naturally to coarsening this cell-level designation as follows: A cell at a coarser level inherits its label based on the labels of the finer-level cells within it. Specifically, a coarser-level cell is labeled as Dirichlet if it encloses any finer-level Dirichlet cell, otherwise as interior if it contains any finer-level interior cell. A cell remains an exterior one if all enclosed finer-level cells are exterior. This discrete domain design enables a straightforward partition of degrees of freedom into two non-overlapping sets for the use of symmetric Vanka and symmetric distributive smoothers. If a cell is proximal to any Dirichlet cell or exterior cell, all the degrees of freedom on this cell are smoothed by symmetric Vanka smoother. Conversely, symmetric distributive smoother is used for all other degrees of freedom.

Every \emph{interior} velocity (on the face between two interior cells) or pressure variable is paired with a discrete equation, according to a finite difference scheme. Dirichlet velocities are similarly paired with a Dirichlet boundary condition. Velocity variables on a face between one interior and one exterior cell are treated as Neumann conditions; we treat those Neumann conditions as zero in our numerical examples, with the understanding that any nonzero conditions can be moved to the right-hand-side of the discretized equations without further alteration of the discrete operator.

\section{Numerical Convergence Experiments}
\label{sec:exp}
In this section, we present the results of our numerical convergence experiments, which compare the performance of our method with the pure multigrid method under several different boundary scenarios. We set the fluid viscosity to $\eta = 10^{-3} m^2/s$, and $\gamma = 10^{-3}$. To enhance the efficiency of Vanka smoothing, we pre-factorize $3^4$ in 2D or $3^6$ in 3D local saddle point problems, based on whether neighboring cells are Dirichlet, exterior or interior. It's straightforward to parallelize the additive Vanka smoother, and an 8-color scheme is applied to parallelize the multiplicative Vanka smoother. All examples are implemented in C++ code, and we utilize Intel Pardiso and Eigen as the direct solvers for the bottom level of the V-cycle and for solving local saddle-point problems in the Vanka smoother respectively. The experiments are performed on a computer with an Intel(R) Xeon(R) CPU E5-2650 v3 @ 2.30GHz and 128 GB of main memory. Paraview is used for postprocessing the results. In the following section, we'll discuss the results in 2D examples and 3D examples separately.

\subsection{2D Numerical Examples}
In this subsection, we explore four 2D examples to assess the performance of our approach. Specifically, we examine the driven cavity example and the Poiseuille flow around a cylinder example to evaluate the convergence rates in comparison to pure multigrid methods. Additionally, the remaining two examples demonstrate the capability of our method to effectively handle more complex scenarios beyond the scope of pure multigrid methods.

\subsubsection{Driven Cavity Example}
\begin{figure}[htb]
\includegraphics[width=\columnwidth]{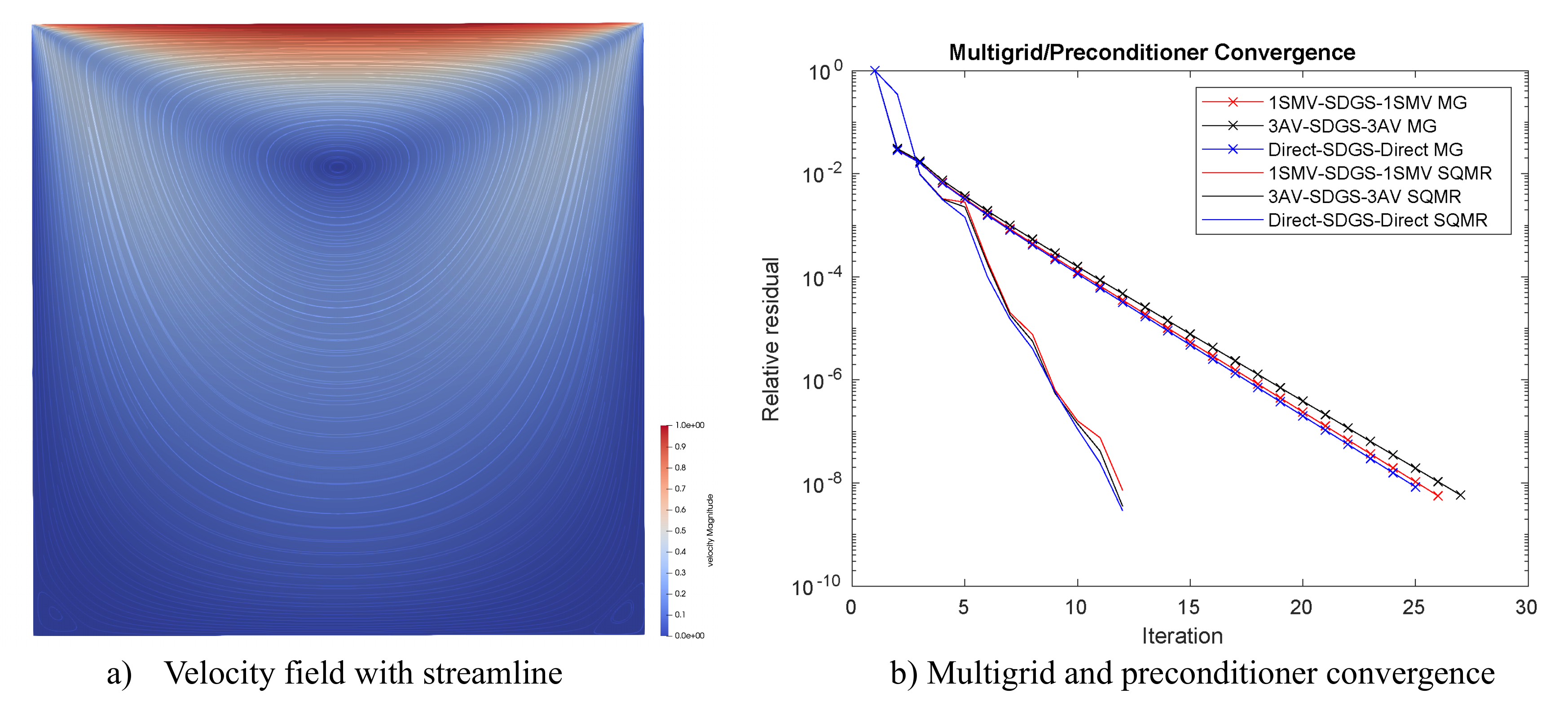}
\centering
\caption{2D $1024^2$ Driven Cavity Example}
\label{fg: 2d-cavity}
\end{figure}

The driven cavity example is a well-established benchmark in the field of computational fluid dynamics. In the 2D scenario, the computational domain is a unit square, where no-slip boundary conditions are imposed on the left, right and bottom sides. On the top side, a zero vertical velocity and a constant horizontal velocity of $u=1.0 m/s$ is imposed. We use a $1024^2$ unit square as our computational domain and an 8-level multigrid as our preconditioner and solver. The solved velocity field with streamline and the convergence of pure multigrid on $\mathbf{\tilde{L}}$ and multigrid-preconditioned SQMR on $\mathbf{L}$ are shown in Fig. \ref{fg: 2d-cavity}. The iterations are terminated when the relative residual goes below $10^{-8}$, and the legend indicates the type and number of smoothers used for the boundary and interior regions. We consistently apply a symmetric distributive Gauss-Seidel smoother only once in each V-cycle to the interior region. This approach significantly reduces the execution time per V-cycle, given that the boundary degrees of freedom are typically less than $3\%$ of the total degrees of freedom. Table. \ref{table: boundary} provides more detailed information about the percentage of the number of boundary degrees of freedom.

The convergence plots reveal that our preconditioner achieves faster convergence compared to the pure multigrid method. Furthermore, the multiplicative Vanka smoother is more efficient than the additive one, even when the latter is applied three times per V-cycle. The convergence rates of both multiplicative and additive methods closely match the ideal scenario where a direct solver is employed for the boundary domains.

\subsubsection{Poiseuille Flow around A Cylinder Example}
\begin{figure}[htb]
\includegraphics[width=\columnwidth]{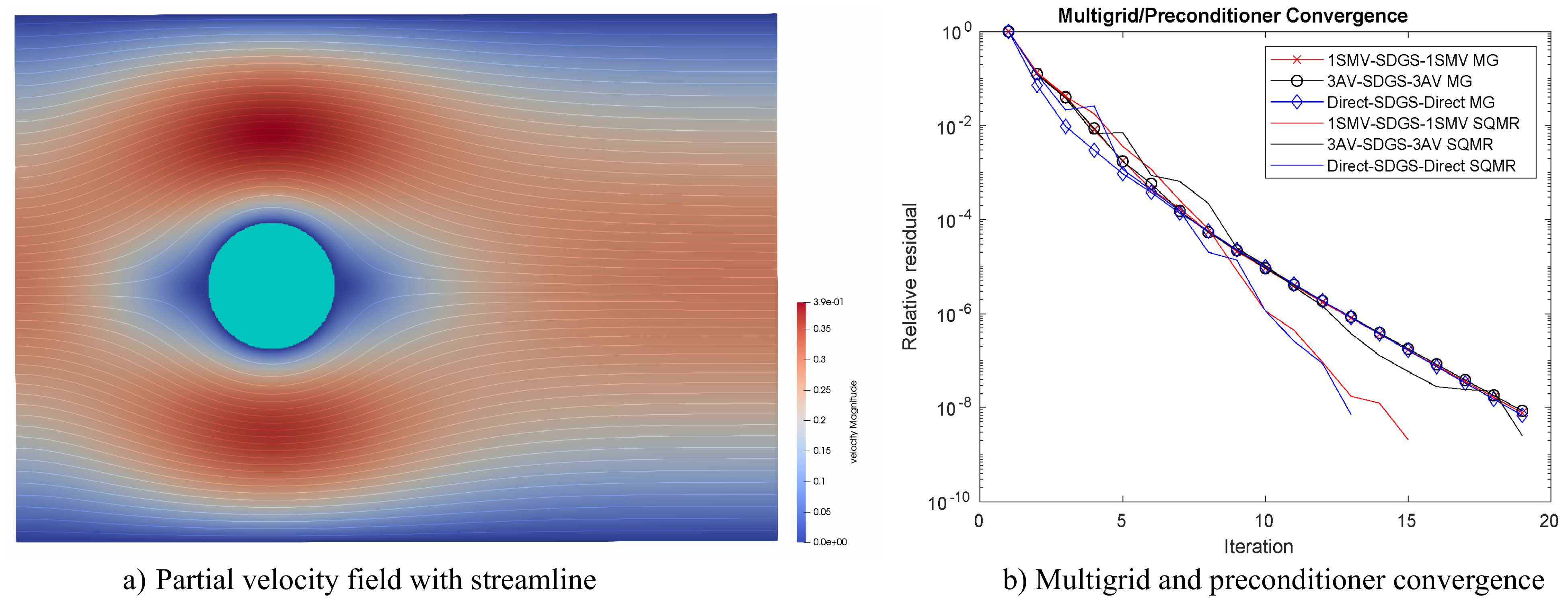}
\centering
\caption{2D $2200 \times 410$ Poiseuille Flow around Cylinder Example}
\label{fg: 2d-cylinder}
\end{figure}

\begin{figure}[htp]
\includegraphics[width=7cm]{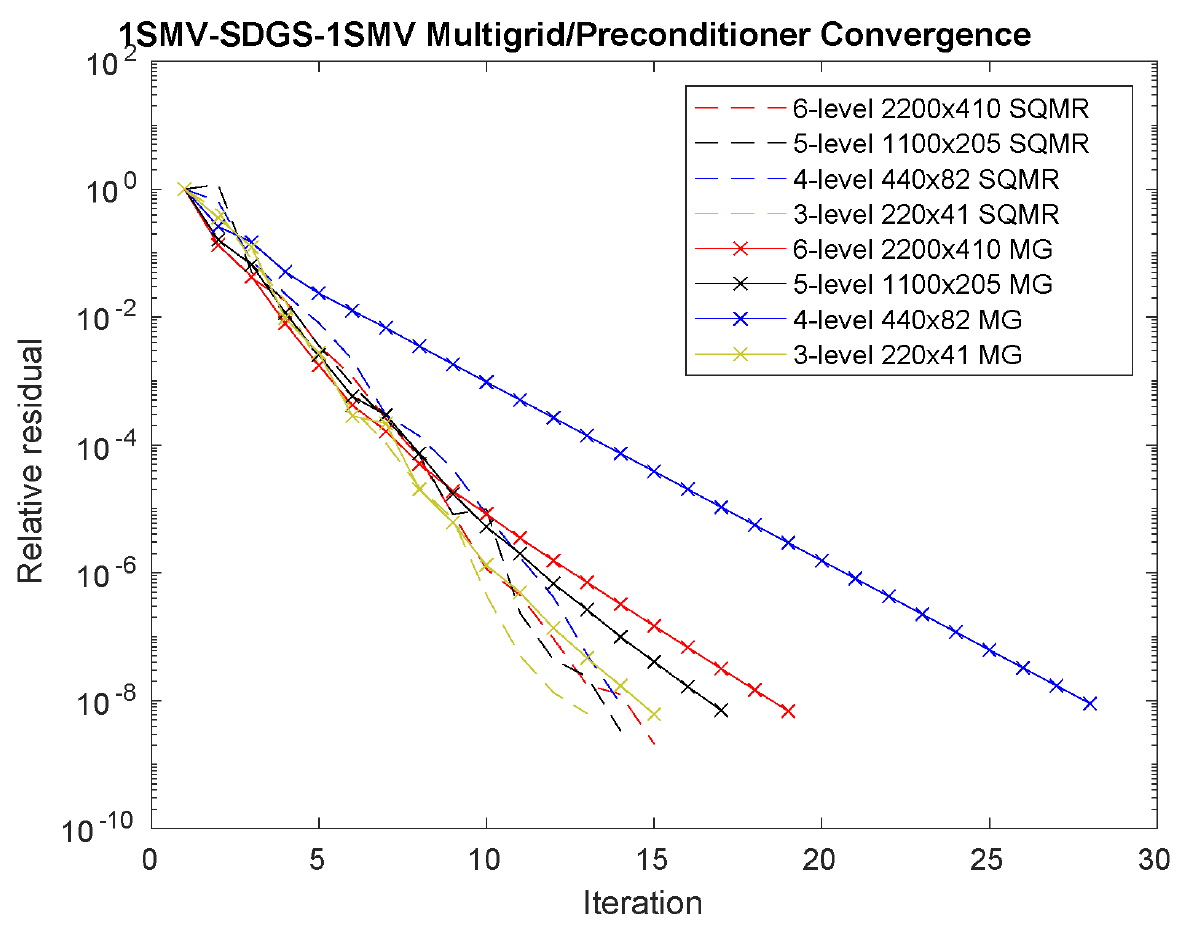}
\centering
\caption{Convergence of 2D Cylinder Example in Different Dimensions}
\label{fg: different-2d-cylinder}
\end{figure}

The case of 2D Poiseuille flow represents the steady, laminar flow of an incompressible fluid through a channel driven by a pressure gradient. This flow scenario is commonly studied in benchmarks \cite{schafer1996benchmark, nicolas2011benchmark}, particularly with the presence of a cylindrical obstacle. In our 2D scenario, the channel dimensions are $2.2m \times 0.41m$, and a cylindrical obstacle is placed at position $(0.2m, 0.2m)$ with a radius of $0.05m$. The inflow boundary conditions are described by
\begin{align}
\begin{split}
    &u(0, y) = \frac{4\bar{u}y(H-y)}{H^2}\\
    &v(0, y) = 0
\end{split}
\end{align}
where $\bar{u}=0.3m/s$ is the average horizontal velocity, and $H=0.41m$ is the width of the channel. On the outflow boundary, a homogeneous Neumann boundary condition is applied, while no-slip boundary conditions are enforced on the top and bottom sides. The computational domain has dimensions $2200 \times 410$, and a 6-level multigrid is used to solve the problem with $h=0.001m$. The partially solved velocity field and the convergence results are shown in Fig. \ref{fg: 2d-cylinder}, which shows that the effectiveness of the symmetric multiplicative Vanka smoother. Additionally, the scalability of our multiplicative preconditioner and multigrid is demonstrated through various domain sizes in Fig. \ref{fg: different-2d-cylinder}.

\subsubsection{Poiseuille Flow around Hollow Square Example}
\begin{figure}[htb]
\includegraphics[width=\columnwidth]{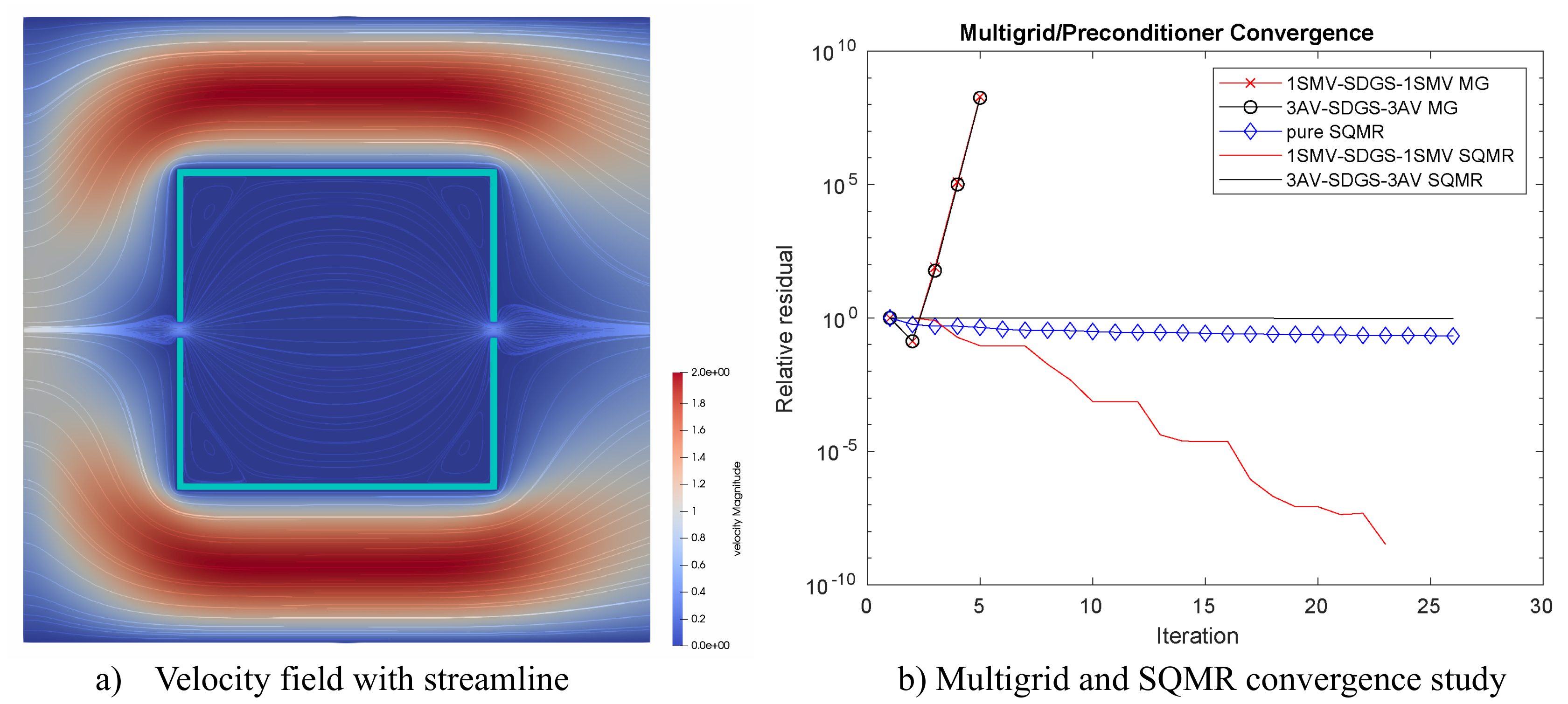}
\centering
\caption{2D $1024^2$ Poiseuille Flow around Hollow Square Example}
\label{fg: 2d-maze}
\end{figure}

In this example, we explore a more complex scenario involving Poiseuille flow around a hollow square obstacle. While pure multigrid suffices for simple applications, it may encounter divergence in more complicated scenarios. Here we use similar boundary conditions as those used in the cylindrical obstacle with $\bar{u}=1.0m/s$. The key difference lies in the presence of a hollow square obstacle, as depicted in Fig. \ref{fg: 2d-maze}. We perform our simulations on a $1024^2$ computational domain with an 8-level multigrid. It is observed that only the multiplicative preconditioner yields fast convergence, while the additive preconditioner and the pure multigrid method fail to converge. This behavior can be attributed to the closing of the hollow square as the multigrid coarsens. Consequently, the solutions derived at coarser levels become increasingly inaccurate. However, when employed as a preconditioner, the multigrid method continues to guide SQMR towards successful convergence.

\subsubsection{Brancher Example}
\begin{figure}[htb]
\includegraphics[width=\columnwidth]{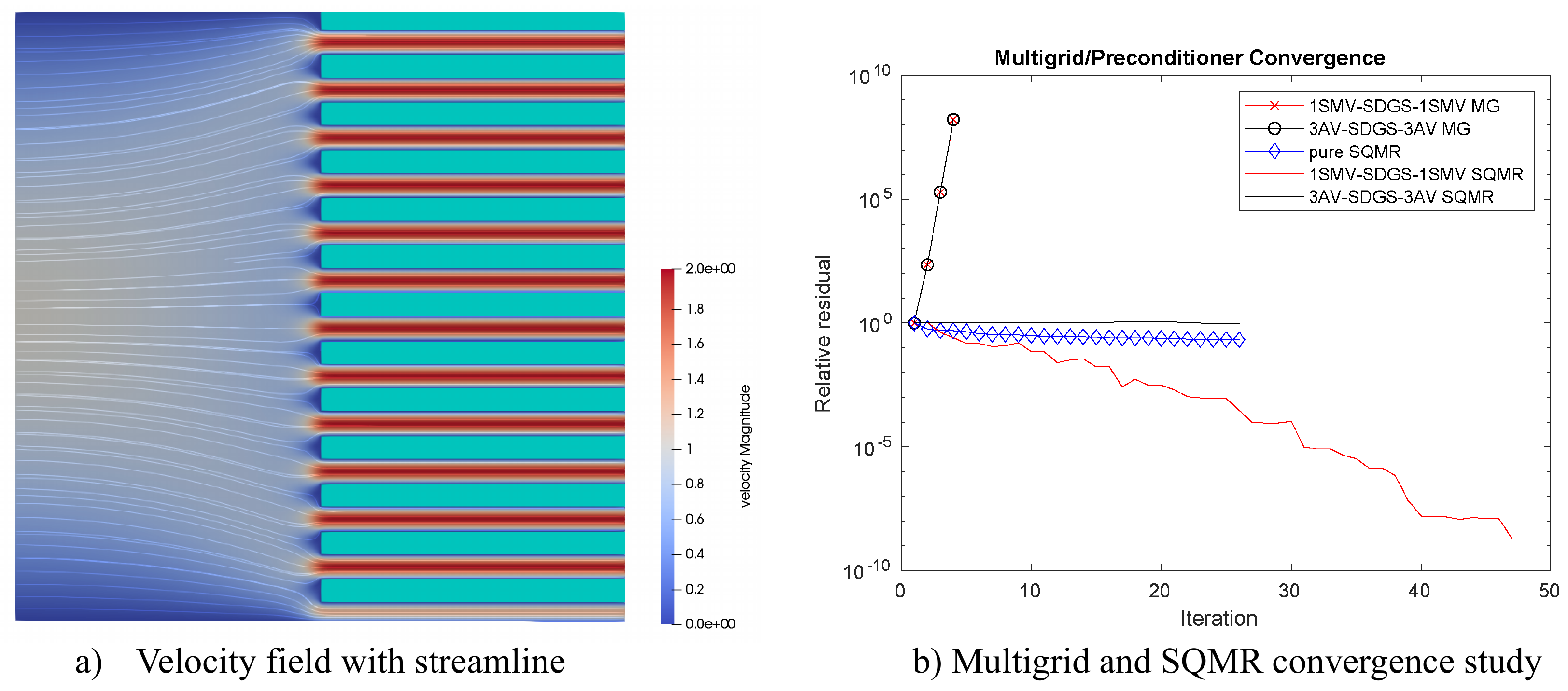}
\centering
\caption{2D $1024^2$ Brancher Example}
\label{fg: 2d-brancher}
\end{figure}

\begin{figure}[H]
\includegraphics[width=7cm]{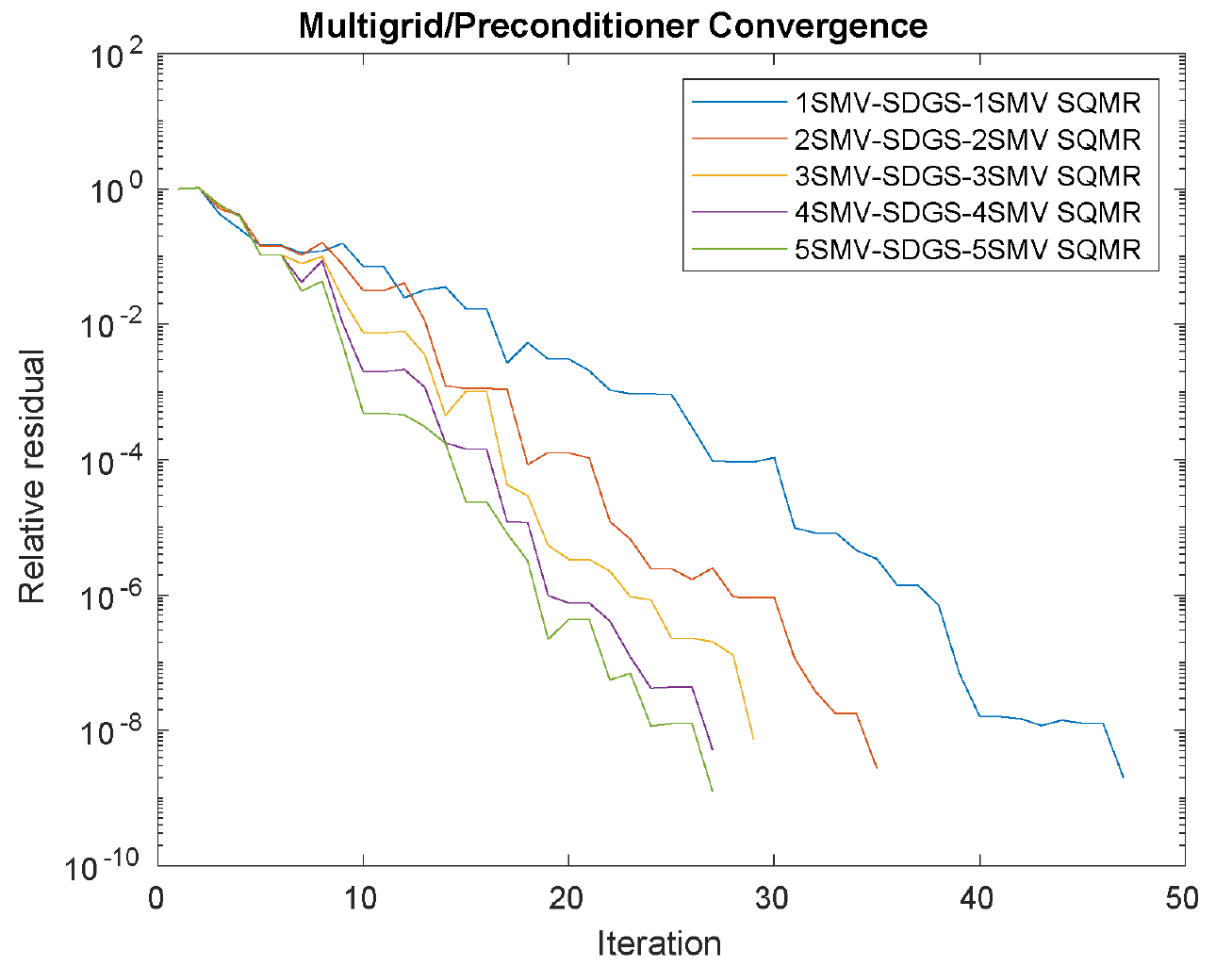}
\centering
\caption{Convergence of 2D Brancher Example with Different Numbers of Boundary Smoothers}
\label{fg: c-2d-brancher-d}
\end{figure}

In the hollow square example, we explored the limitations of pure multigrid methods. Now we present a more practical example to emphasize the defects of pure multigrid methods and the reliability of the multiplicative preconditioner. We keep the same boundary conditions and computational domain as those used in the hollow square obstacle. Instead of a single obstacle, we apply multiple rectangular obstacles near the outflow boundary to branch a single inflow into several outflows, forming a structure referred to as a ``brancher". As shown in Fig. \ref{fg: 2d-brancher}, our multiplicative preconditioner converges very fast. In contrast, the additive preconditioner and the pure multigrid methods still struggle to converge. The hollow square example and the brancher example clearly show the effectiveness and reliability of the multiplicative preconditioner. This fact leads us to focus on the multiplicative methods in our upcoming 3D examples, without considering the additive method. Additionally, we investigate the impact of the number of boundary smoother iterations on the convergence rate of the multiplicative preconditioner in Fig. \ref{fg: c-2d-brancher-d}. The result indicates that increasing the numbers of boundary smoother iterations leads to faster convergence rates.

\subsection{3D Numerical Examples}
In this subsection, we consider four 3D examples using only multiplicative preconditioner and multigrid. The driven cavity example and the Poiseuille flow with a cylindrical obstacle example show that our multigrid preconditioner achieves comparable convergence rates to pure multigrid methods. The last two examples highlight the capability of our approach to address more complicated scenarios beyond the capabilities of pure multigrid solvers.

\subsubsection{Driven Cavity Example}
\begin{figure}[htb]
\includegraphics[width=\columnwidth]{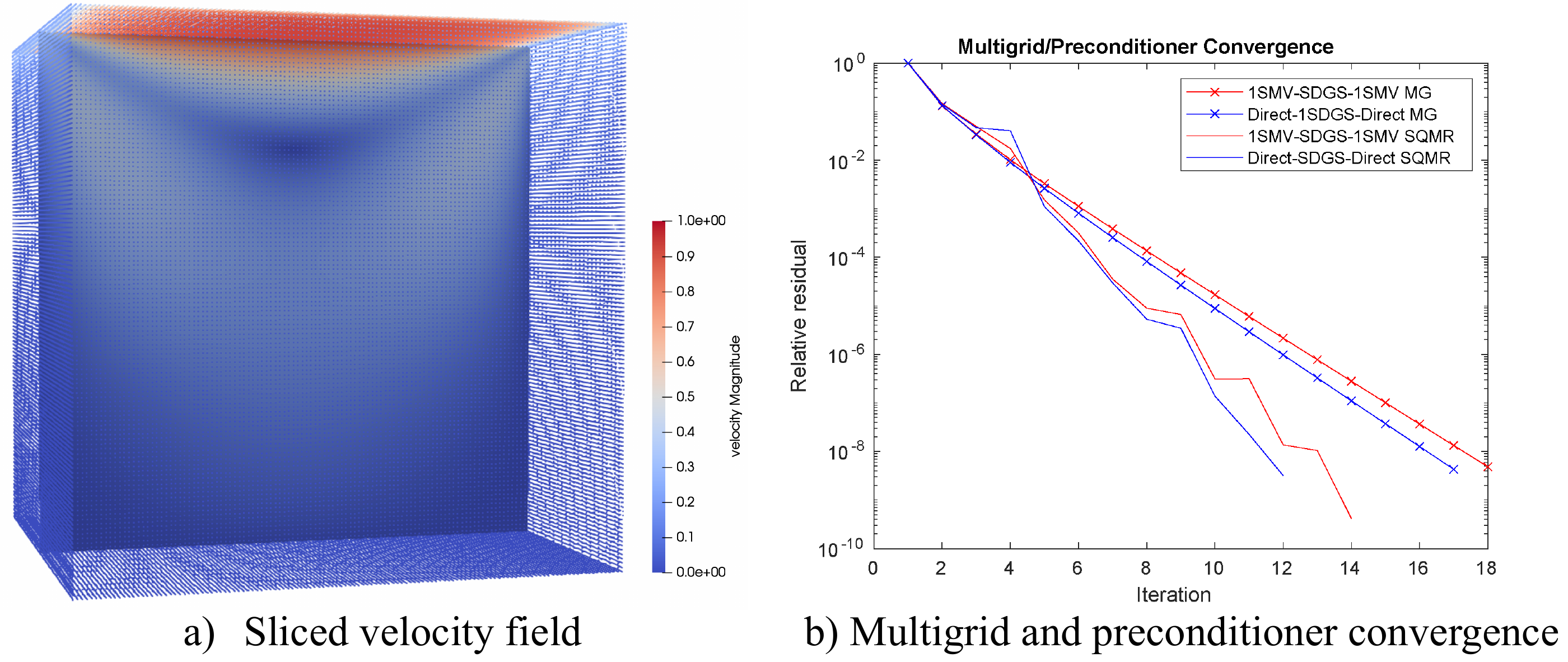}
\centering
\caption{3D $128^3$ Driven Cavity Example}
\label{fg: 3d-cavity}
\end{figure}

The 3D driven cavity example is similar to the 2D version. In this case, we consider a unit cube instead of a unit square as our computational domain, applying no-slip boundary conditions to all sides except the top. The top side is subject to a constant horizontal velocity of $u=1.0m/s$. A 5-level multiplicative multigrid is used within a $128^3$ unit cube. We present the sliced velocity field and the convergence results in Fig. \ref{fg: 3d-cavity}. We can see that our preconditioner method converges faster than the multigrid method, closely resembling where the boundary domains are solved using a direct solver.

\subsubsection{Poiseuille Flow around Cylinder Example}
\begin{figure}[htb]
\includegraphics[width=\columnwidth]{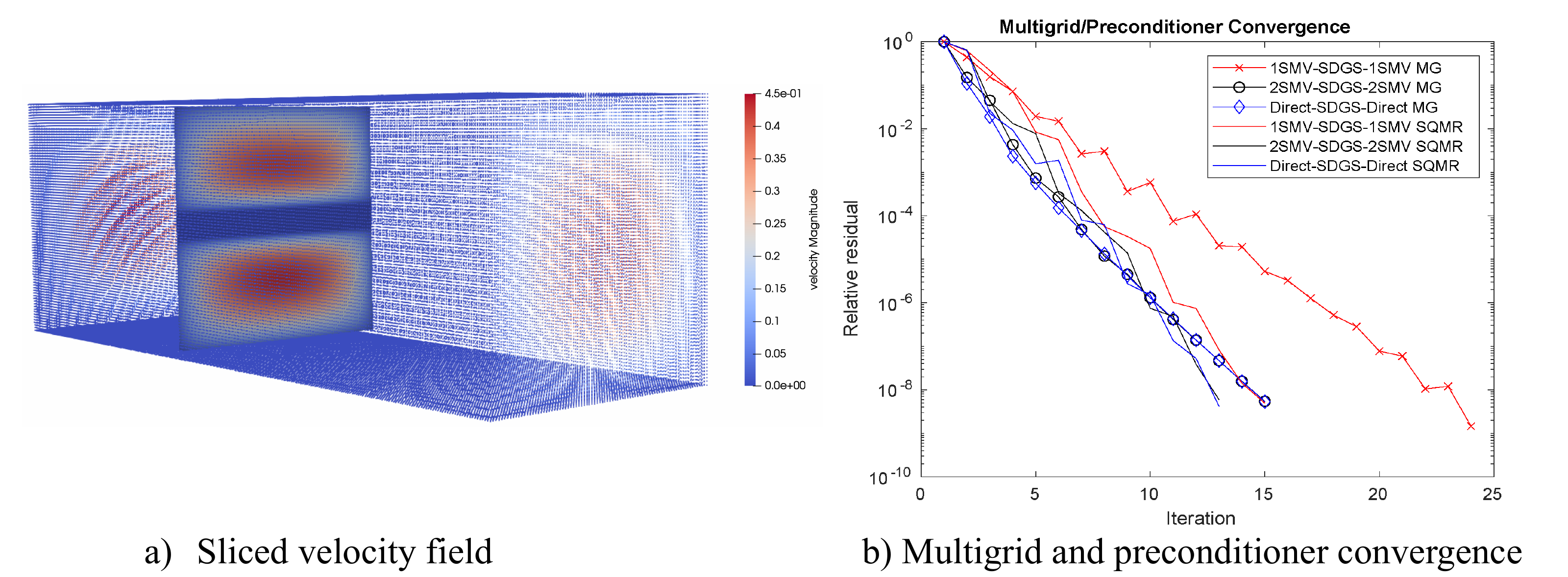}
\centering
\caption{3D $255 \times 82 \times 82$ Cylinder Example}
\label{fg: 3d-cylinder}
\end{figure}

In the 3D scenario, we also consider Poiseuille flow around a cylinder. The channel dimensions are $1.275m \times 0.41m \times 0.41m$, and a cylindrical obstacle is positioned along the $(0.5m, 0.2m, z)$ axis with a radius of $0.05m$. The inflow boundary conditions are described as
\begin{align}
\begin{split}
    &u(0, y, z) = \frac{16\bar{u}yz(H_y - y)(H_z - z)}{H^4}\\
    &v(0, y, z) = 0\\
    &w(0, y, z) = 0
\end{split}
\end{align}
where $\bar{u}=0.45m/s$ represents the average horizontal velocity, and $H_y=H_z=0.41m$ are the inflow widths of the channel. Homogeneous Neumann boundary conditions are applied on the outflow boundary, while the surrounding sides are subject to no-slip conditions. The computational domain has dimensions of $255 \times 82 \times 82$, and a 4-level multigrid is used to solve the problem with $h=0.005m$. Fig. \ref{fg: 3d-cylinder} b) shows the sliced velocity field and the convergence results. In this example, applying our symmetric multiplicative Vanka smoother only once does not yield convergence rates that match the ideal scenario where a direct solver is utilized on the boundary domains. Instead, we achieve the ideal convergence by increasing the number of smoothers applied to the boundary to 2.

\subsubsection{Brancher Example}
\begin{figure}[htb]
\includegraphics[width=14cm]{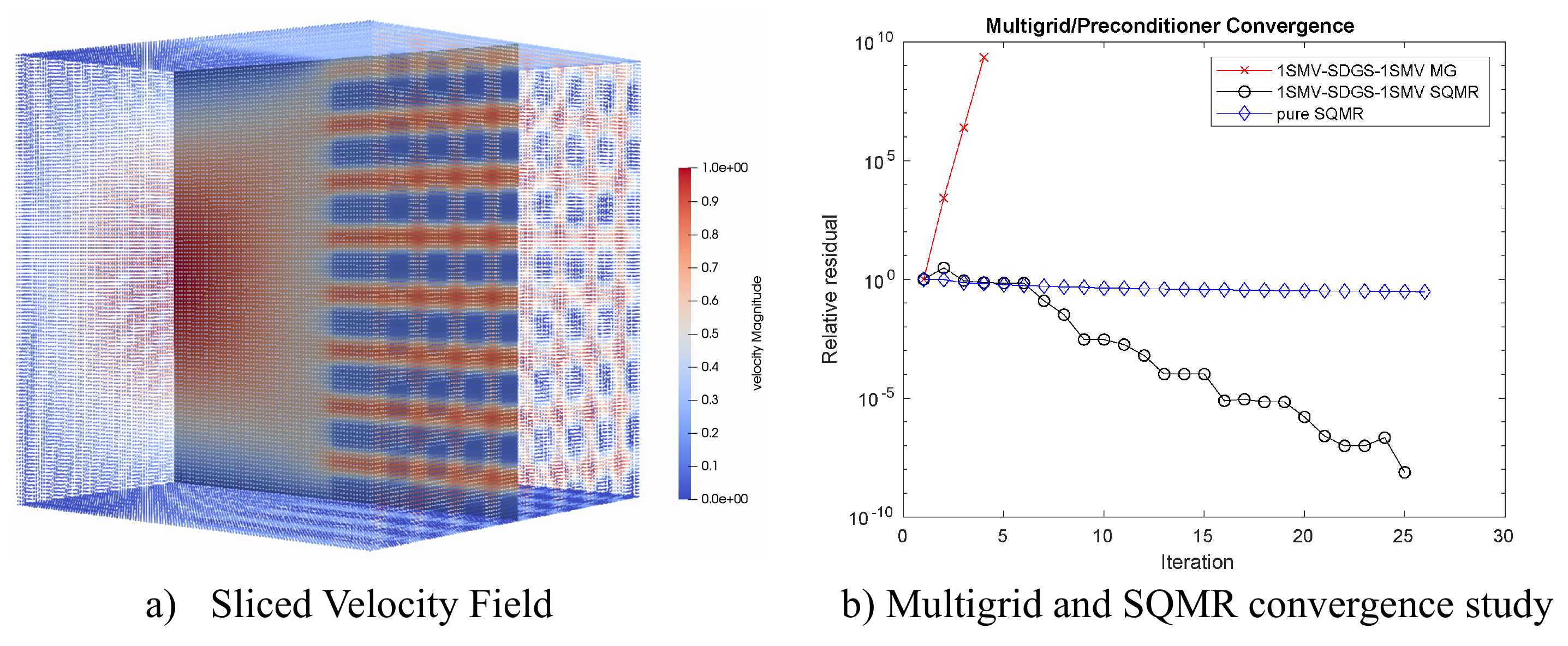}
\centering
\caption{3D $128^3$ Brancher Example}
\label{fg: 3d-brancher}
\end{figure}

To highlight the effectiveness of our preconditioner, we present a 3D brancher example. We keep the same boundary condition as 3D Poiseuille flow around cylinder example except that we substitute the cylinder with multiple cuboids near the outflow boundary and use $\bar{u}=1.0m/s$. A 5-level multigrid is used within a $128^3$ unit cube. The sliced velocity field and the convergence results are shown in Fig. \ref{fg: 3d-brancher}. We can see that our multiplicative preconditioner converges very fast but the pure SQMR and multigrid method fails. This example shows the advantage of our approach in 3D case.

\subsubsection{Porous Example}
\begin{figure}[htb]
\includegraphics[width=\columnwidth]{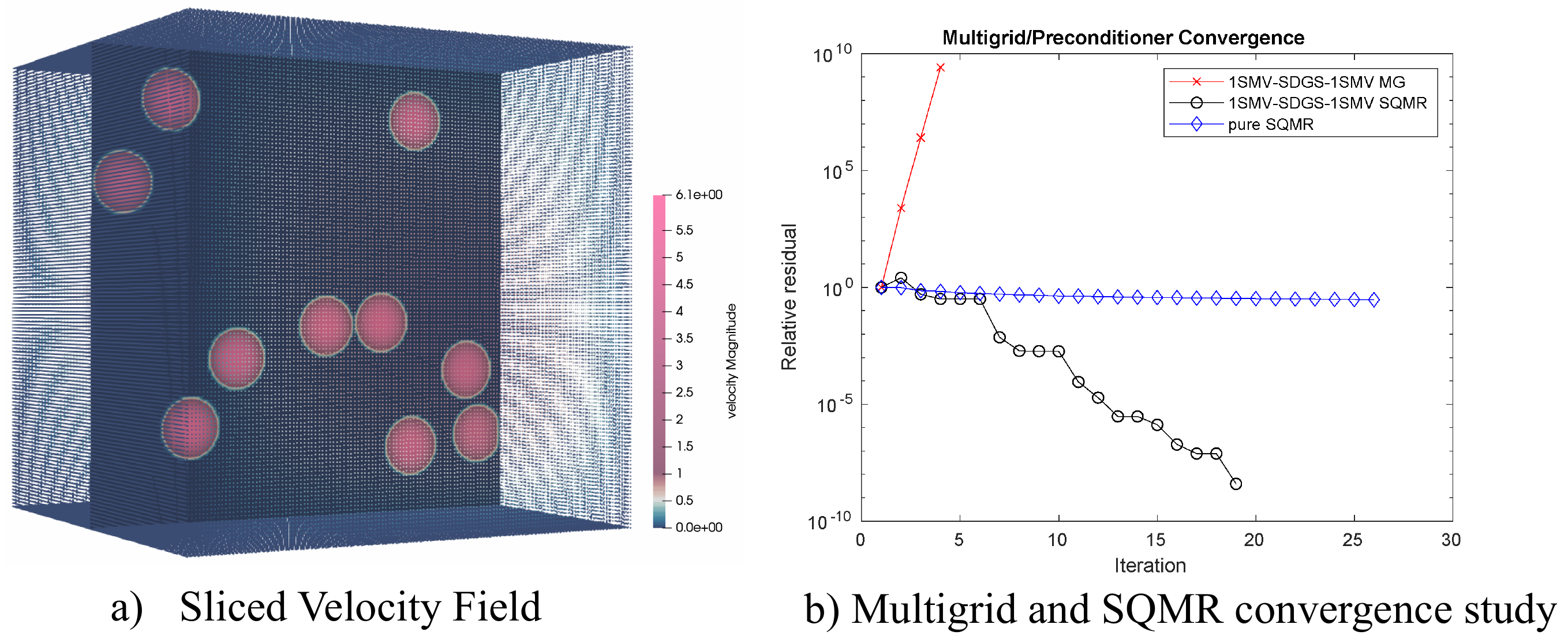}
\centering
\caption{3D $128^3$ Brancher Example}
\label{fg: 3d-porous}
\end{figure}

The final example we show is the porous example. In this example, we use the same boundary condition and computational domain as brancher example, but we use a porous layer as shown in Fig. \ref{fg: 3d-porous} a) to replace the multiple cuboid obstacles. our multiplicative preconditioner exhibits rapid convergence, while the pure SQMR and multigrid methods prove ineffective, as illustrated in Fig. \ref{fg: 3d-porous} b). This example further exemplifies the robustness and efficiency of our multiplicative preconditioner.

\begin{table}[htb]
\resizebox{\columnwidth}{!}{%
\begin{tabular}{
>{\columncolor[HTML]{DAE8FC}}c ccc}
\cellcolor[HTML]{C0C0C0}\textbf{2D examples} & \cellcolor[HTML]{C0C0C0}\textbf{\#total degrees of freedom}                        & \cellcolor[HTML]{C0C0C0}\textbf{\#boundary degrees of freedom}                     & \cellcolor[HTML]{C0C0C0}\textbf{boundary percentage}                     \\
$1024^2$ Driven Cavity                       & 3143680                                                                            & 12272                                                                              & 0.38\%                                                                   \\
$2200 \times 410$ Cylinder                   & 2680020                                                                            & 18493                                                                              & 0.69\%                                                                   \\
$1100 \times 205$ Cylinder                   & 669372                                                                             & 9243                                                                               & 1.38\%                                                                   \\
$440 \times 82$ Cylinder                     & 106812                                                                             & 3693                                                                               & 3.45\%                                                                   \\
$220 \times 41$ Cylinder                     & 26580                                                                              & 1843                                                                               & 6.93\%                                                                   \\
$1024^2$ Hollow Square                       & 3463768                                                                            & 29305                                                                              & 0.85\%                                                                   \\
$1024^2$ Brancher                            & 2546560                                                                            & 61077                                                                              & 2.39\%                                                                   \\
\cellcolor[HTML]{C0C0C0}\textbf{3D examples} & \cellcolor[HTML]{C0C0C0}{\color[HTML]{333333} \textbf{\#total degrees of freedom}} & \multicolumn{1}{l}{\cellcolor[HTML]{C0C0C0}\textbf{\#boundary degrees of freedom}} & \multicolumn{1}{l}{\cellcolor[HTML]{C0C0C0}\textbf{boundary percentage}} \\
$128^3$ Driven Cavity                        & 8339456                                                                            & 385580                                                                             & 4.62\%                                                                   \\
$255 \times 82 \times 82$ Cylinder           & 6788860                                                                            & 425856                                                                             & 6.27\%                                                                   \\
$128^3$ Brancher                             & 8013824                                                                            & 1119980                                                                            & 13.98\%                                                                  \\
$128^3$ Porous                               & 9156798                                                                            & 614304                                                                             & 6.71\%                                                                  
\end{tabular}%
}
\caption{The numbers of total and boundary degrees of freedom}
\label{table: boundary}
\end{table}

\section{Stability and Efficiency Analysis}
\label{sec:additional_exp}
We present some additional experiments to assess the stability of our method under problems of different resolutions and conditioning, and multigrid cycle schemes.

\subsection{Resolution and Conditioning}
\begin{figure}[htb]
\centering
\includegraphics[width=7cm]{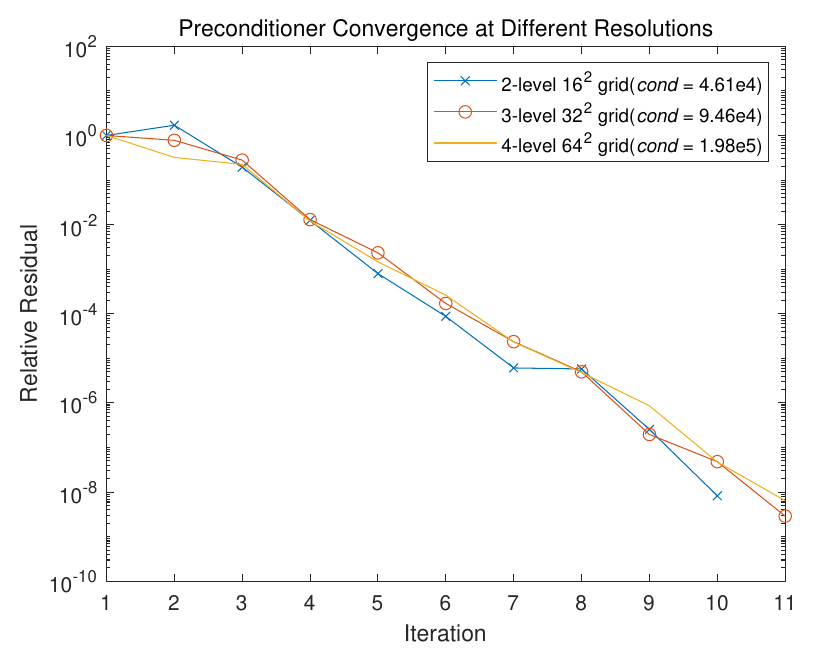}
\caption{Convergence of 2D Driven Cavity Example at Different Resolutions}
\label{fg: condition_number}
\end{figure}

We seek to investigate the effect of domain resolution on the efficacy of our scheme. Although the increased spatial resolution does impact the condition number of the modified system matrix $\mathbf{\tilde{L}}$ used in multigrid, a healthy multi-resolution scheme would be expected to be consistently convergent, largely independent of resolution, even if the condition number itself increases. We illustrate the convergence of the 2D driven cavity example at different resolutions as shown in Fig. \ref{fg: condition_number}. We observe consistent convergence rates in this refinement analysis. We should indicate that the \emph{geometric} complexity of a domain is another factor that can certainly affect convergence, especially in the presence of fine geometric features that are not well resolved or consistently represented across resolutions.

\subsection{Cycle Schemes}
\begin{figure}[htb]
\centering
\includegraphics[width=7cm]{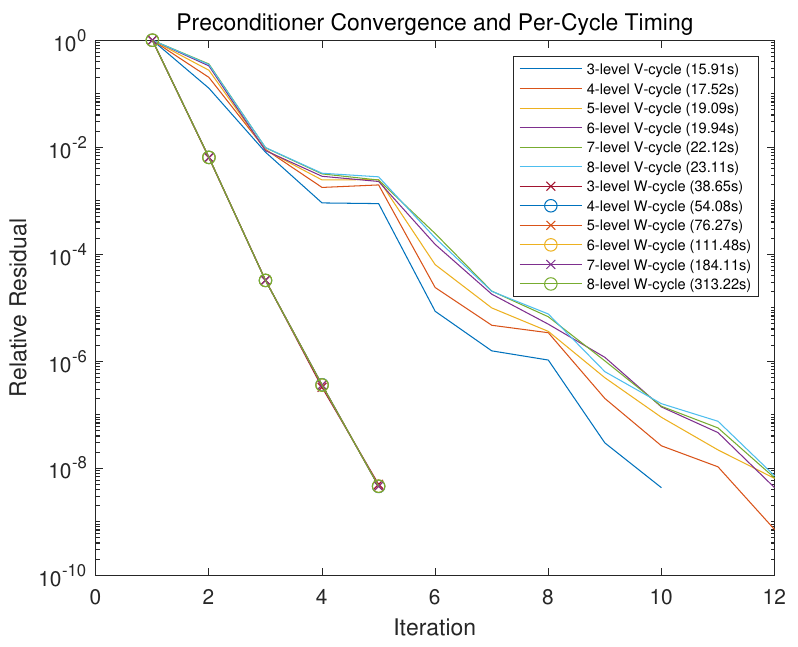}
\caption{Different Cycle Schemes on 2D $1024^2$ Driven Cavity Example}
\label{fg: wcycle_cavity}
\end{figure}

\begin{figure}[htb]
\centering
\includegraphics[width=7cm]{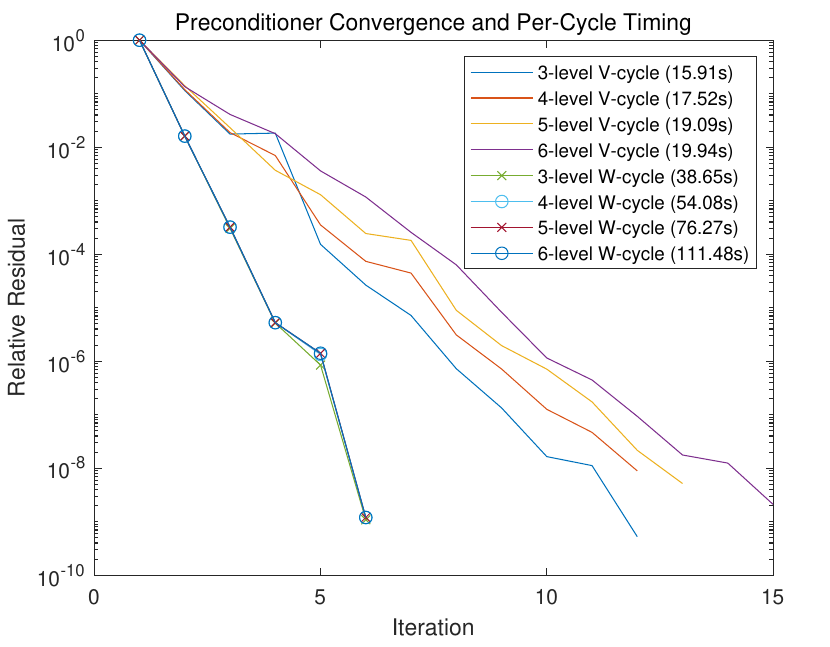}
\caption{Different Cycle Schemes on 2D $2200 \times 410$ Poiseuille Flow around A Cylinder Example}
\label{fg: wcycle_cylinder}
\end{figure}

We also design experiments using different cycle schemes. Specifically, we compare the convergence rates by using W-cycle and V-cycle on $1024^2$ Driven Cavity Example and $2200 \times 410$ Poiseuille Flow around Cylinder Example as shown in Fig. \ref{fg: wcycle_cavity} and \ref{fg: wcycle_cylinder}. The results indicate that the W-cycle can indeed accelerate convergence rates (independent of runtime) in the cases where even pure multigrid methods converge. However, the increased time required for the W-cycle does not fully compensate for the improved convergence rates. Additionally, the deep-level W-cycle scheme fails to converge in our $1024^2$ Poiseuille flow around a hollow square example and brancher example. This suggests that the V-cycle is a more reliable scheme and potentially economical in terms of runtime.


\section{Conclusions}
\label{sec:conclusion}
In this work, we propose and evaluate a novel multigrid preconditioner for solving the Stokes equations in both 2D and 3D scenarios. Our approach leverages the strengths of a symmetric distributive smoother and Vanka smoother, which are well-suited for handling interior and boundary regions respectively. By combining these smoothers carefully, we achieve an efficient symmetric multigrid preconditioner which can be used by the SQMR method. We demonstrate the efficacy of our scheme in numerical experiments in both 2D and 3D scenarios. Our method achieves comparable convergence rates to pure multigrid methods, and also maintains fast convergence even when pure multigrid methods struggle with divergence. The application of our preconditioner as a robust and reliable solver, particularly in challenging 3D cases, shows its potential to significantly enhance the efficiency of solving the Stokes equations. Our future endeavors will focus on GPU acceleration to enhance the scalability of our methods, allowing for the efficient resolution of larger-scale problems.

\appendix

\section{Symmetry proof for symmetric distributive Gauss-Seidel}
\label{sec:pDGS}
\begin{lemma}
Consider a sequence $\mathbf{x}^{(k+1)}=\mathbf{A}\mathbf{x}^{(k)}+\mathbf{b}$ with $\mathbf{x}^{(0)} = \mathbf{0}$. The general formula is $\mathbf{x}^{(n)} = \sum_{i=0}^{n-1}\mathbf{A}^i\mathbf{b}$.
\label{lemma1}
\end{lemma}

This conclusion is obvious by using mathematical induction. \qedsymbol

\begin{theorem}
In symmetric DGS smoother, when $\mathbf{L}$ is symmetric and $\mathbf{x}^{(0)}=\mathbf{0}$, then $\mathbf{x}^{(n)}=\mathbf{C}\mathbf{x}^{(0)}$ and $\mathbf{C}$ is a symmetric matrix.
\label{theorem1}
\end{theorem}

Splitting $\mathbf{L}=\mathbf{G}+\mathbf{H}$ with a lower triangular matrix and a strictly upper triangular matrix, GS can be written as
\begin{equation}
    \mathbf{x}^{(k+1)}=\mathbf{x}^{(k)}+\mathbf{G}^{-1}(\mathbf{b}-\mathbf{L}\mathbf{x}^{(k)})
\end{equation}
After introducing the distributive preconditioner $\mathbf{M}$ and defining $\mathbf{\tilde{G}}$ and $\mathbf{\tilde{H}}$ as the lower and strictly upper triangular part of $\mathbf{L}\mathbf{M}$ instead of $\mathbf{L}$, the distributive GS can be written as
\begin{equation}
    \mathbf{x}^{(k+1)}=\mathbf{x}^{(k)}+\mathbf{M}\mathbf{\tilde{G}}^{-1}(\mathbf{b}-\mathbf{L}\mathbf{x}^{(k)})  
\end{equation}
For the symmetric DGS smoother, the smoothing step can be written as
\begin{equation}
    \begin{split}
    &\mathbf{x}^{(k+\frac{1}{2})}=\mathbf{x}^{(k)}+\mathbf{M}\mathbf{\tilde{G}}^{-1}(\mathbf{b}-\mathbf{L}\mathbf{x}^{(k)})\\
    &\mathbf{x}^{(k+1)}=\mathbf{x}^{(k+\frac{1}{2})}+\mathbf{\tilde{G}}^{-T}\mathbf{M}^{T}(\mathbf{b}-\mathbf{L}\mathbf{x}^{(k+\frac{1}{2})})
    \end{split}
    \label{eq:SDGS}
\end{equation}
Substitute the first equation in Equation \ref{eq:SDGS} into the second one and denote $\mathbf{P}=\mathbf{M}\mathbf{\tilde{G}}^{-1}+\mathbf{\tilde{G}}^{-T}\mathbf{M}^T-\mathbf{\tilde{G}}^{-T}\mathbf{M}^T\mathbf{L}\mathbf{M}\mathbf{\tilde{G}}^{-1}$, and we get
\begin{equation}
    \mathbf{x}^{(k+1)} = (\mathbf{I} - \mathbf{P}\mathbf{L})\mathbf{x}^{(k)} + \mathbf{P}\mathbf{b}
\end{equation}
By using Lemma \ref{lemma1}, the general formula of $\mathbf{x}^{(n)}$ is
\begin{equation}
    \mathbf{x}^{(n)} = \sum_{i=0}^{n-1}(\mathbf{I}-\mathbf{P}\mathbf{L})^i\mathbf{P}\mathbf{b}=\sum_{i=0}^{n-1}\sum_{j=0}^i{i \choose j}(-\mathbf{P}\mathbf{L})^i\mathbf{P}\mathbf{b}=\mathbf{C}\mathbf{b}
\end{equation}
Since $(-\mathbf{P}\mathbf{L})^i\mathbf{P}$ is always symmetric when $\mathbf{L}$ and $\mathbf{P}$ are symmetric, $\mathbf{C}$ is symmetric. \qedsymbol

\section{Symmetry proof for additive Vanka smoother}
\label{sec:pAV}
\begin{theorem}
In additive Vanka smoother, when $\mathbf{L}$ is symmetric and $\mathbf{x}^{(0)}=\mathbf{0}$, $\mathbf{x}^{(n)}=\mathbf{C}\mathbf{x}^{(0)}$ and $\mathbf{C}$ is a symmetric matrix.
\label{theorem2}
\end{theorem}
We know that the series of additive Vanka smoother can be written as
\begin{equation}
    \label{eq:av-relaxation}
    \mathbf{x}^{(k+1)}=\mathbf{x}^{(k)}+\mathbf{S}_{AV}(\mathbf{b}-\mathbf{L}\mathbf{x}^{(k)}) = (\mathbf{I} - \mathbf{S}_{AV}\mathbf{L})\mathbf{x}^{(k)} + \mathbf{S}_{AV}\mathbf{b}
\end{equation}
By using Lemma \ref{lemma1}, the general formula of $\mathbf{x}^{(n)}$ is
\begin{equation}
    \mathbf{x}^{(n)} = \sum_{i=0}^{n-1}(\mathbf{I}-\mathbf{S}_{AV}\mathbf{L})^i\mathbf{S}_{AV}\mathbf{b}=\sum_{i=0}^{n-1}\sum_{j=0}^i{i \choose j}(-\mathbf{S}_{AV}\mathbf{L})^i\mathbf{S}_{AV}\mathbf{b}=\mathbf{C}\mathbf{b}
\end{equation}
Since $(-\mathbf{S}_{AV}\mathbf{L})^i\mathbf{S}_{AV}$ is always symmetric when $\mathbf{L}$ and $\mathbf{S}_{AV}$ are symmetric, $\mathbf{C}$ is symmetric. \qedsymbol

\section{Symmetry proof for symmetric multiplicative Vanka smoother}
\label{sec:pMV}
\begin{theorem}
In symmetric multiplicative Vanka smoother, when $\mathbf{L}$ is symmetric and $\mathbf{x}^{(0)}=\mathbf{0}$, then $\mathbf{x}^{(n)}=\mathbf{C}\mathbf{x}^{(0)}$ and $\mathbf{C}$ is a symmetric matrix.
\label{theorem3}
\end{theorem}

We know that the series of multiplicative Vanka smoother can be written as
\begin{equation}
    \label{eq:mv-relaxation}
    \mathbf{x}^{(k+1)}=\mathbf{x}^{(k)}+\mathbf{S}_{MV}(\mathbf{b}-\mathbf{L}\mathbf{x}^{(k)}) = (\mathbf{I} - \mathbf{S}_{MV}\mathbf{L})\mathbf{x}^{(k)} + \mathbf{S}_{MV}\mathbf{b}
\end{equation}
For the symmetric multiplicative Vanka smoother, the smoothing step can be written as
\begin{equation}
    \begin{split}
    &\mathbf{x}^{(k+\frac{1}{2})}=(\mathbf{I} - \mathbf{S}_{MV}^{(1)}\mathbf{L})\mathbf{x}^{(k)} + \mathbf{S}_{MV}^{(1)}\mathbf{b}\\
    &\mathbf{x}^{(k+1)}=(\mathbf{I} - \mathbf{S}_{MV}^{(2)}\mathbf{L})\mathbf{x}^{(k)} + \mathbf{S}_{MV}^{(2)}\mathbf{b}
    \end{split}
    \label{eq:SMV}
\end{equation}
where $\mathbf{S}_{MV}^{(1)}$ is solving in the order from domain $D_1$ to $D_n$, and $\mathbf{S}_{MV}^{(2)}$ is solving in the order from domain $D_n$ to $D_1$.
Before we prove that this smoother is symmetric, let's prove that $\mathbf{S}_{MV}^{(2)}=(\mathbf{S}_{MV}^{(1)})^T$. Denoting that $\mathbf{K}_i=\mathbf{C}_i^T(\mathbf{C}_i\mathbf{L}\mathbf{C}_i^T)^{-1}\mathbf{C}_i$, We can rewrite the transpose of $\mathbf{S}_{MV}^{(1)}$ as
\begin{equation}
    \begin{split}
    (\mathbf{S}_{MV}^{(1)})^T &= ([\mathbf{I}-\prod_{i = 1}^{n}(\mathbf{I} - \mathbf{K}_i\mathbf{L})]\mathbf{L}^{-1})^T=\mathbf{L}^{-1}[\mathbf{I}-\prod_{i = 1}^{n}(\mathbf{I} - \mathbf{K}_i\mathbf{L})]^{T}\\
    &=\mathbf{L}^{-1}-\mathbf{L}^{-1}\prod_{i = n}^{1}(\mathbf{I} - \mathbf{L}\mathbf{K}_i)=\mathbf{L}^{-1}-[\prod_{i = n}^{1}(\mathbf{I} - \mathbf{K}_i\mathbf{L})]\mathbf{L}^{-1}\\
    &=[\mathbf{I}-\prod_{i = n}^{1}(\mathbf{I} - \mathbf{K}_i\mathbf{L})]\mathbf{L}^{-1}=\mathbf{S}_{MV}^{(2)}
    \end{split}
\end{equation}
Substitute the first equation in Equation \ref{eq:SMV} into the second one and denote $\mathbf{P}=\mathbf{S}_{MV}^{(1)}+[\mathbf{S}_{MV}^{(1)}]^T-[\mathbf{S}_{MV}^{(1)}]^T\mathbf{L}\mathbf{S}_{MV}^{(1)}$, we can get
\begin{equation}
    \label{C4}
    \mathbf{x}^{(k+1)} = (\mathbf{I} - \mathbf{P}\mathbf{L})\mathbf{x}^{(k)} + \mathbf{P}\mathbf{b}
\end{equation}
By using Lemma \ref{lemma1}, the general formula of $\mathbf{x}^{(n)}$ is
\begin{equation}
    \mathbf{x}^{(n)} = \sum_{i=0}^{n-1}(\mathbf{I}-\mathbf{P}\mathbf{L})^i\mathbf{P}\mathbf{b}=\sum_{i=0}^{n-1}\sum_{j=0}^i{i \choose j}(-\mathbf{P}\mathbf{L})^i\mathbf{P}\mathbf{b}=\mathbf{C}\mathbf{b}
\end{equation}
Since $(-\mathbf{P}\mathbf{L})^i\mathbf{P}$ is always symmetric when $\mathbf{L}$ and $\mathbf{P}$ are symmetric, $\mathbf{C}$ is symmetric. \qedsymbol

\section{Symmetry proof for symmetric integration of smoother}
\label{sec:pVDV}
\begin{lemma}
Consider a sequence $\mathbf{x}^{(k+1)}=\mathbf{A}\mathbf{x}^{(k)}+\mathbf{b}$. The general formula is $\mathbf{x}^{(n)} = \mathbf{A}^n\mathbf{x}^{(0)} + \sum_{i=0}^{n-1}\mathbf{A}^i\mathbf{b}$.
\label{lemma2}
\end{lemma}

It's obvious by using mathematical induction. \qedsymbol

\begin{theorem}
For smoother shown in Alg. \ref{alg:smoothing}, when $\mathbf{L}$ is symmetric and $\mathbf{x}^{(0)}=\mathbf{0}$, then $\mathbf{x}^{(n)}=\mathbf{C}\mathbf{x}^{(0)}$ and $\mathbf{C}$ is a symmetric matrix.
\label{theorem4}
\end{theorem}

Since we decouple the unknowns into two sets $V_1$ and $V_2$, we can reorder our system equations into
\begin{equation}
    \mathbf{L}\mathbf{x} =
    \begin{pmatrix}
        \mathbf{L}_{11} & \mathbf{L}_{12}\\
        \mathbf{L}_{21} & \mathbf{L}_{22}\\
    \end{pmatrix}
    \begin{pmatrix}
        \mathbf{x}_1\\
        \mathbf{x}_2\\
    \end{pmatrix}
    = \begin{pmatrix}
        \mathbf{b}_1\\
        \mathbf{b}_2\\
    \end{pmatrix}
    =\mathbf{b}
\end{equation}
where $x_1$ and $x_2$ represent unknowns for $V_1$ and $V_2$ separately.

For symmetric Vanka smoothing, it can be seen as solving the system
\begin{equation}
    \mathbf{L}_{11}\mathbf{x}_1 = \mathbf{b}_1 - \mathbf{L}_{12}\mathbf{x}_2
\end{equation}
And for symmetric distributive smoothing, it can be seen as solving the system
\begin{equation}
    \mathbf{L}_{22}\mathbf{x}_2 = \mathbf{b}_2 - \mathbf{L}_{21}\mathbf{x}_1
\end{equation}

For the first symmetric Vanka smoothing, $\mathbf{x}^{(0)}=\mathbf{0}$. By Theorem \ref{theorem2}
or \ref{theorem3}, we know that $[\mathbf{x}_1, \mathbf{x}_2]^T=[\mathbf{C}_1\mathbf{b}_1, \mathbf{0}]^T$ after smoothing, where $\mathbf{C}_1$ has the form $\sum_{i=0}^{n-1}\sum_{j=0}^i{i \choose j}(-\mathbf{P}\mathbf{L})^i\mathbf{P}\mathbf{b}$ with some symmetric matrix $\mathbf{P}$.

For the symmetric Vanka smoothing, $\mathbf{x}^{(0)}=[\mathbf{C}_1\mathbf{b}_1, \mathbf{0}]^T$. By Theorem \ref{theorem1}, we know that $[\mathbf{x}_2, \mathbf{x}_2]^T=[\mathbf{C}_1\mathbf{b}_1, \mathbf{C}_2(\mathbf{b}_2 - \mathbf{L}_{21}\mathbf{C}_1\mathbf{b}_1)]^T$ after smoothing, where $\mathbf{C}_2$ is also some symmetric matrix.

For the second symmetric Vanka smoothing, $\mathbf{x}^{(0)}=[\mathbf{C}_1\mathbf{b}_1, \mathbf{C}_2(\mathbf{b}_2 - \mathbf{L}_{21}\mathbf{C}_1\mathbf{b}_1)]^T$. By Lemma \ref{lemma2} and Equation \ref{C4}, we know that 
\begin{equation}
    \begin{pmatrix}
        \mathbf{x}_1\\
        \mathbf{x}_2\\
    \end{pmatrix}
    = \begin{pmatrix}
        \mathbf{C}_1(\mathbf{b}_1 - \mathbf{L}_{12}\mathbf{C}_2(\mathbf{b}_2 - \mathbf{L}_{21}\mathbf{C}_1\mathbf{b}_1)) + (\mathbf{I}-\mathbf{P}\mathbf{L})^n\mathbf{C}_1\mathbf{b}_1\\
        \mathbf{C}_2(\mathbf{b}_2 - \mathbf{L}_{21}\mathbf{C}_1\mathbf{b}_1)
    \end{pmatrix}
\end{equation}
after smoothing. And we can rewrite it in the matrix form
\begin{equation}
    \begin{split}
    \begin{pmatrix}
        \mathbf{x}_1\\
        \mathbf{x}_2\\
    \end{pmatrix}
    &= \begin{pmatrix}
        \mathbf{C}_1 + \mathbf{C}_1\mathbf{L}_{12}\mathbf{C}_2\mathbf{L}_{21}\mathbf{C}_1+(\mathbf{I}-\mathbf{P}\mathbf{L})^n\mathbf{C}_1 & -\mathbf{C}_1\mathbf{L}_{12}\mathbf{C}_2\\
        -\mathbf{C}_2\mathbf{L}_{21}\mathbf{C}_1 & \mathbf{C}_{2}\\
    \end{pmatrix}
    \begin{pmatrix}
        \mathbf{b}_1\\
        \mathbf{b}_2\\
    \end{pmatrix}\\
    &= \mathbf{C}\mathbf{b}\\
    \end{split}
\end{equation}
Since $\mathbf{C}_1$ has the form $\sum_{i=0}^{n-1}\sum_{j=0}^i{i \choose j}(-\mathbf{P}\mathbf{L})^i\mathbf{P}\mathbf{b}$ with some symmetric matrix $\mathbf{P}$, $(\mathbf{I}-\mathbf{P}\mathbf{L})^n\mathbf{C}_1$ is still symmetric. Therefore, $\mathbf{C}$ is symmetric. \qedsymbol\\

\noindent
\textbf{Acknowledgements}\\
This work was supported in part by NSF Grant IIS-2106768.\\

 \bibliographystyle{elsarticle-num} 
 \bibliography{MGPSQMR}





\end{document}